\documentclass[12pt]{article}

\usepackage{amssymb, amsmath, amsthm, bbm}
\usepackage{cite, enumerate}
\usepackage[arrow, matrix, curve]{xy}
\usepackage{xcolor}
\usepackage{graphicx}
\usepackage{tikz}

\newcommand{\N}{\mathbb{N}}
\newcommand{\R}{\mathbb{R}}
\newcommand{\E}{\mathbb{E}}
\newcommand{\K}{\mathcal{K}}

\newcommand{\SP}{\mathbb{S}}
\renewcommand{\S}{\mathbb{S}}

\DeclareMathOperator{\BD}{bd}
\DeclareMathOperator{\ORTHO}{O}

\DeclareMathOperator{\INT}{int}

\newcounter{abschnitt}
\newtheorem{satz}{Theorem}

\newtheorem{theorem}{Theorem}[abschnitt]
\newtheorem{koro}[theorem]{Corollary}
\newtheorem{prop}[theorem]{Proposition}
\newtheorem{lem}[theorem]{Lemma}

\newtheorem*{defi}{Definition}

\renewenvironment{quote}{\list{}{\leftmargin=0.62in\rightmargin=0.62in}\item[]}{\endlist}

\newcounter{saveeqn}
\newcommand{\alpheqn}{\setcounter{saveeqn}{\value{abschnitt}}
\renewcommand{\theequation}{\mbox{\arabic{saveeqn}.\arabic{equation}}}}
\newcommand{\reseteqn}{\setcounter{equation}{0}
\renewcommand{\theequation}{\arabic{equation}}}

\begin{document}

\phantom{a}

\vspace{-1.7cm}

\begin{center}
\begin{Large} {\bf Spherical Centroid Bodies} \\[0.4cm] \end{Large}

\begin{large} Florian Besau, \,\! Thomas Hack, \\[0.2cm] Peter Pivovarov, \,\! Franz E. Schuster\end{large}
\end{center}

\vspace{-0.8cm}

\begin{quote}
\footnotesize{ \vskip 1cm \noindent {\bf Abstract.} The spherical centroid body of a centrally-symmetric convex body in the Euclidean unit sphere is introduced.
Two alternative definitions -- one geometric, the other probabilistic in nature -- are given and shown to lead to the same objects. The geometric approach is then used to
establish a number of basic properties of spherical centroid bodies, while the probabilistic approach inspires the proof of a spherical analogue of the classical polar Busemann--Petty centroid inequality.}
\end{quote}

\vspace{0.5cm}

\centerline{\large{\bf{ \setcounter{abschnitt}{1} \arabic{abschnitt}. Introduction}}}

\alpheqn

\vspace{0.5cm}

For an origin-symmetric convex body $K$ (that is, a compact, convex set with nonempty interior) in an $n$-dimensional linear vector space,
the centroids of the intersections of $K$ with half-spaces through the origin form the surface of its centroid body $\Gamma K$. In the case $n = 3$, this construction
first explicitly appeared in a paper by Blaschke \textbf{\cite{Blaschke:1917}}, where he conjectured that the ratio of the volume of a body to that
of its centroid body attains its maximum for ellipsoids. This conjecture was confirmed by Petty \textbf{\cite{Petty:1961}} (who also coined the name centroid bodies), by
reinterpreting Busemann's random simplex inequality \textbf{\cite{Busemann:1953}} as what would become known  as the \emph{Busemann--Petty centroid inequality}.
Combining it with the famous \linebreak Blaschke--Santal\'o inequality, leads to an isoperimetric inequality for the product of the volumes of a body and its \emph{polar} centroid body, called the \emph{polar Busemann--Petty \linebreak centroid inequality}. An $L_p$ analogue of the latter inequality was established by Lutwak and Zhang \textbf{\cite{LZ1997}}, who showed that it provides an $L_p$ extension of the Blaschke--Santal\'o inequality. A stronger $L_p$ Busemann--Petty centroid inequality was obtained later in \textbf{\cite{LYZ2000a}} and \textbf{\cite{Campi:Gronchi02a}}
(see also, \textbf{\cite{habschu09, Ivaki2016, Lutwak:1986, Lutwak:1990, LYZ:2010, Nguyen:2018, zhu2012}} for further generalizations).

The definition of centroid bodies implies that the operator $\Gamma$ commutes with all non-degenerate linear transformations. Therefore, both the Busemann--Petty centroid inequality and its polar version belong to the affine geometry of convex bodies, where, in particular, they were very recently used to establish a series of affine invariant Sobolev-type inequalities (see
\textbf{\cite{DeNapoli:etal:2018, HJM:2016, HJM:2018, HJM:2019, Nguyen:2016}}). In fact, such affine isoperimetric and analytic inequalities have proven to be significantly more powerful than their Euclidean counterparts (see e.g., \textbf{\cite{habschu19, LYZ2002, zhang99}}). But it is not just the impact of these isoperimetric inequalities that would make it difficult to overstate the importance of the notion of centroid bodies; it is also the fact that they have naturally appeared in many different contexts ranging from asymptotic geometric analysis (see e.g., \textbf{\cite{Brazitikos:etal:2014, Klartag:Milman:2012, Milman:Pajor:1989, Paouris:2006}}), geometric tomography
(see e.g., \textbf{\cite{Gardner:2006, Yaskin:Yaskina:2006, Ivaki2017}}), and integral geometry (see e.g., \textbf{\cite{Haberl:2012, Ludwig:2007}}) to recently even Finsler geometry (see \textbf{\cite{Bernig:2014}}) and information theory (see \textbf{\cite{Paouris:Werner:2012}}).

\pagebreak

This article belongs to a line of research of recent origins, dealing with the question of which \emph{affine} constructions and isoperimetric inequalities from linear vector spaces allow for generalizations to other spaces of, say, constant curvature (then no longer compatible with the affine group but rather the isometry group of the respective space). The starting point for these investigations was the proof of a spherical counterpart of one of the best known affine isoperimetric inequalities, the Blaschke--Santal\'o inequality, by Gao, Hug, and Schneider \textbf{\cite{Gao:Hug:Schneider:2003}}. Since this seminal paper appeared, considerable efforts have been invested to obtain further results in the same spirit
(see \textbf{\cite{Besau:Schuster:2016, Besau:Werner:2016, Besau:Werner:2018, Besau:Ludwig:Werner:2018, Dann:etal:2018, Oliker:2007, Yaskin:2006}}). However, despite these substantial inroads, this area of research is still in its infancy with several important threads left open.

Here, we attack one of these threads by introducing a spherical analogue of the centroid body of a centrally-symmetric convex body in the Euclidean unit sphere. For our (first) definition of spherical centroid bodies, we mimic Blaschke's geometric approach to centroid bodies in linear vector spaces laid out at the beginning of this introduction. Combining our geometric definition with the gnomonic projection, naturally leads to centroid bodies (in the tangent linear space) with respect to a specific weight. These \emph{weighted centroid bodies} will allow us to deduce several basic properties of spherical centroid bodies such as continuity in the Hausdorff metric and injectivity as well as the fact that, like in the linear setting, all spherical centroid bodies are $C^2$-smooth and have everywhere positive Gau\ss\ curvature.

Weighted centroid bodies and their associated isoperimetric inequalities have only recently become critical tools in high-dimensional probability, where they are used to establish important small-ball probabilities for marginals of probability measures (see \textbf{\cite{Dann:etal:2016, Paouris:2012a, Paouris:Valettas:2014}}). Sharp volume inequalities for weighted centroid bodies and their polars were first established in \textbf{\cite{Paouris:2012}} and \textbf{\cite{Cordero:etal:2015, Paouris:Pivovarov:2012}}, the latter of which are based on random approximations of these bodies (see Lemma \ref{4:approxGammaMu} for details). The obtained results for empirical analogues of centroid bodies yield stronger stochastic dominance inequalities and are part of a research program going back to \textbf{\cite{Paouris:Pivovarov:2012}} aimed at proving random extensions of classical inequalities from convexity such as the Brunn-Minkowski inequality and its relatives (see \textbf{\cite{Paouris:Pivovarov:2018}} for more information).

With our first main result we show that, as in the linear setting, there is an alternative probabilistic approach towards spherical centroid bodies. Hereby, they are obtained from the spherical convex hull of centroids of uniformly distributed random points in the given convex body by letting the number of points tend to infinity. Similarly, by considering \emph{randomized} weighted centroid bodies, we are able to establish our second main result, a \emph{sharp spherical polar Busemann--Petty centroid inequality}.

\vspace{0.2cm}

In order to state our results more precisely, let $\mathbb{S}^n$ denote the unit $n$-sphere of Euclidean space $\mathbb{R}^{n+1}$, where we always assume that $n \geq 2$.
For $u \in \mathbb{S}^n$, we write $\mathbb{S}_u^+ = \{v \in \mathbb{S}^n: u \cdot v \geq 0$\} for the \emph{closed} hemisphere centered around $u$ and $\mathbb{S}_u:=\mathbb{S}^n \cap u^{\bot}$ for its boundary.

A set $A \subseteq \mathbb{S}^n$ is called (\emph{spherically}) \emph{convex}, if $\mathrm{rad}\,A := \{\lambda u: \lambda \geq 0, u \in A\}$ is convex in $\mathbb{R}^{n+1}$.
We call $K \subseteq \mathbb{S}^n$ a (\emph{spherically}) \emph{convex body} if $K \neq \mathbb{S}^n$ is closed, convex, and has nonempty interior (with respect to $\mathbb{S}^n$) and we say $K$ is \emph{proper} if $K$ is contained in some open hemisphere. Finally, for $A \subseteq \mathbb{S}^n$, we denote its interior and boundary relative to $\mathbb{S}^n$ by $\mathrm{int}\,A$ and $\mathrm{bd}\,A$, respectively, and its \emph{spherical centroid} by $c_s(A)$ whenever it exists (see Section 2 for details).

\begin{defi}
For a convex body $K \subseteq \mathbb{S}^n$ which is centrally-symmetric with center $e \in \mathbb{S}^n$, we define its spherical centroid body $\Gamma_s K$ by
\[\mathrm{bd}\,\Gamma_s K := \{c_s(K\cap\mathbb{S}_u^+):u\in\mathbb{S}_e\}.\]
\end{defi}

In Section 3 we show that $\Gamma_s K$ is indeed a well defined proper spherically convex body which is centrally-symmetric with the same center as $K$. Our main tool to see this, is the well known \emph{gnomonic projection}. In order to recall its definition, we denote by $\mathbb{R}^n_u$ the linear subspace in $\mathbb{R}^{n+1}$ orthogonal to $u \in \mathbb{S}^n$. Then the gnomonic projection (with respect to $u$), is defined by
\[g_u: \mathrm{int}\,\mathbb{S}^+_u \rightarrow \mathbb{R}^n_u, \qquad g_u(v)= \frac{v}{u \cdot v}-u.  \]
It is particularly useful in spherical geometry since it maps proper convex bodies in $\mathrm{int}\,\mathbb{S}^+_u$ to convex bodies in $\mathbb{R}^n_u$. More importantly for our purposes, we will show that spherical centroid bodies are mapped by the gnomonic projection to centroid bodies with respect to the following weight function,
\[\psi: \mathbb{R}^n \rightarrow \mathbb{R}^+, \qquad \psi(x) = \left (1+\|x\|^2\right)^{-\frac{n+2}{2}}.  \]
Here, $\|\cdot\|$ denotes the standard Euclidean norm in $\mathbb{R}^n$. In order to state the definition of weighted centroid bodies (in linear vector spaces), let $H_u^+ = \{x \in \mathbb{R}^n: u \cdot x \geq 0\}$ denote the closed halfspace in $\mathbb{R}^n$ with exterior normal $u \in \mathbb{S}^{n-1}$. Then, for an even finite Borel measure $\mu$ on $\mathbb{R}^n$ with positive density and an origin-symmetric convex body $L \subseteq \mathbb{R}^n$, the \emph{$\mu$-centroid body} $\Gamma_{\mu}L$ of $L$ can be defined by
\begin{equation} \label{introdefmycent}
\mathrm{bd}\,\Gamma_{\mu} L := \{c_{\mu}(L\cap H_u^+):u\in\mathbb{S}^{n-1}\}.
\end{equation}
Here, we write $c_{\mu}$ for the center of mass with respect to $\mu$. In Section 3, we then show that for a proper convex body $K \subseteq \mathbb{S}^n$ centered around $e \in \mathbb{S}^n$, the following critical relation holds,
\begin{equation} \label{gnomGammas}
g_e(\Gamma_s K) = \Gamma_{\widehat{\tau}}\, g_e(K),
\end{equation}
where $\widehat{\tau}$ is the absolutely continuous measure on $\mathbb{R}^n_e$ (w.r.t.\ Lebesgue measure) with density $\psi$.
While relation (\ref{gnomGammas}) already allows us to establish a number of basic properties of spherical centroid bodies, we also consider them under a different, probabilistic perspective to justify our definition on the one hand, and to gain a deeper understanding of this new notion on the other hand.

In the linear setting, the probabilistic approach towards centroid bodies was first noted in \textbf{\cite{Paouris:Pivovarov:2012}}, and can be described as follows. Given an origin-symmetric convex body $L \subseteq \mathbb{R}^n$ and $N \in \mathbb{N}$ independent random points $X_1,\dotsc, X_N$ distributed uniformly in $L$, define the (random) convex body
\begin{equation} \label{randcent}
\Gamma(X_1,\ldots, X_N) := \frac{1}{N}\sum_{i=1}^N [-X_i,X_i] = \mathrm{conv}\, \{c(\pm X_1,\ldots,\pm X_N)\},
\end{equation}
where $[-X_i,X_i]$ denotes the line segment joining $\pm X_i$, the sum is the standard Minkowski addition, and $c(x_1,\ldots,x_N)$ denotes the usual centroid of finitely many points in $\mathbb{R}^n$.
The crucial observation from \textbf{\cite{Paouris:Pivovarov:2012}} is that $\Gamma(X_1,\ldots, X_N)$ converges almost surely in the Hausdorff metric to the centroid body $\Gamma L$ as $N$ tends to infinity.

Although there is no natural analogue of Minkowski addition on $\mathbb{S}^n$, both the convex hull and centroids of finite point sets do have natural counterparts. In order to mimic definition (\ref{randcent}) on $\mathbb{S}^n$, we can therefore use the second equation in (\ref{randcent}), but here we replace $- v$ with the geodesic reflection of $v$ about a point $e \in \mathbb{S}^n$, that is, $v \mapsto v^e := -v + 2(v \cdot e)e$, and use the abbreviation $v^{(e)}$ for $\{v,v^e\}$.

\begin{defi} For a proper finite set $\{u_1,\ldots, u_N\} \subseteq \mathbb{S}^n$ contained in $\mathrm{int}\,\mathbb{S}_e^+$ for some $e \in \mathbb{S}^n$, we define
\[ \Gamma_{s,e}(u_1,\ldots,u_N) := \mathrm{conv}\,\left \{c_s\left (u_1^{(e)},\ldots, u_N^{(e)}\right)\right\}.\]
\end{defi}

Our first main result is a spherical version of the random approximation result for centroid bodies from \textbf{\cite{Paouris:Pivovarov:2012}}.

\begin{satz} \label{1:thmApprox}
Let $K \subseteq \mathbb{S}^n$ be a spherically convex body which is centrally-symmetric with center $e \in \mathbb{S}^n$.
If $U_1, \ldots, U_N$ are independent random unit vectors distributed uniformly in $K$, then
\[\Gamma_{s,e}(U_1,\ldots,U_N) \to \Gamma_s K\]
almost surely in the spherical Hausdorff metric as $N$ tends to infinity.
\end{satz}

In order to prove an isoperimetric inequality for spherical centroid bodies w.r.t.\ spherical Lebesgue measure $\sigma$, we can use (\ref{gnomGammas}) and consider the pushforward measure $\widehat{\sigma}$ of $\sigma$ under gnomonic projection of the weighted centroid bodies $\Gamma_{\widehat{\tau}}\,L$.
By refining ideas and techniques from \textbf{\cite{Cordero:etal:2015}}, we are able to show that for any origin-symmetric convex body $L \subseteq \mathbb{R}^n$,
\begin{equation} \label{sigmaGammatau}
\widehat{\sigma} \left (\Gamma_{\widehat{\tau}}^\circ\, L\right) \leq \widehat{\sigma}\left (\Gamma_{\widehat{\tau}}^\circ\, B_L^{\widehat{\tau}}\right),
\end{equation}
where $\Gamma_{\widehat{\tau}}^\circ\, L$ denotes the usual polar body of $\Gamma_{\widehat{\tau}}\, L$ in $\mathbb{R}^n$ (see Section 2) and $B_L^{\widehat{\tau}}$ is the
Euclidean ball in $\mathbb{R}^n$ centered at the origin such that $\widehat{\tau}(B_L^{\widehat{\tau}}) = \widehat{\tau}(L)$.

\pagebreak

By pulling (\ref{sigmaGammatau}) back to $\mathbb{S}^n$ with the inverse gnomonic projection, we obtain our second main result, a \emph{spherical polar Busemann--Petty centroid inequality}. In order to state it, let $\tau$ be the absolutely continuous measure on $\S^n$ (w.r.t.\ to $\sigma$) with density $ u \mapsto |e\cdot u|$ (then, $\widehat{\tau}$ is the pushforward of $\tau$ under gnomonic projection). Moreover, for a spherically convex body $K \subseteq \mathbb{S}^n$ centered around $e \in \mathbb{S}^n$, we write $K^* = \{u\in \mathbb{S}^n: u\cdot v \leq 0 \mbox{ for all } v\in K\}$ for the spherical polar of $K$ and $C_K^\tau$ for the spherical cap centered at $e$ such that $\tau(K)=\tau(C_K^\tau)$.

\begin{satz} \label{1:thmIneq}
If $K \subseteq \mathbb{S}^n$ is a spherically convex body which is centrally-symmetric with center $e \in \mathbb{S}^n$, then
\[\sigma(\Gamma_s^* K) \leq \sigma(\Gamma_s^* C_K^\tau).\]
\end{satz}

\vspace{1cm}

\centerline{\large{\bf{ \setcounter{abschnitt}{2}
\arabic{abschnitt}. Background material}}}

\reseteqn \alpheqn \setcounter{theorem}{0}

\vspace{0.6cm}

In the following we first recall basic definitions and facts from spherical geometry. However, the main part of this section is devoted to the gnomonic projection and spherical centroids as well as their interplay for which we prove several auxiliary results needed in the next sections. As a general reference for this section we recommend \textbf{\cite{Besau:Schuster:2016}} and \textbf{\cite{Glasauer:1996}}.

The usual \emph{spherical distance} between two points on the $n$-dimensional Euclidean unit sphere $\mathbb{S}^n$ is given by $d_s(u, v) = \arccos(u \cdot v)$, $u, v \in \mathbb{S}^n$.
For $r > 0$ and $A \subseteq \mathbb{S}^n$, we write $A_r$ for the set of all points with distance at most $r$ from $A$. In particular, we denote by $C_r(u):=\{u\}_r$ the spherical cap of radius $r\geq 0$ centered at $u \in \mathbb{S}^n$. The Hausdorff distance between closed sets $A, B \subseteq \mathbb{S}^n$ is given by
\[\delta_s(A,B) = \min\{ 0 \leq r \leq \pi: A \subseteq B_r \mbox{ and } B \subseteq A_r\}.\]

We write $\mathcal{K}(\mathbb{S}^n)$ for the space of \emph{spherically convex bodies} endowed with the Hausdorff metric. As usual, the convex hull of $A \subseteq \mathbb{S}^n$ is the intersection of all convex sets in $\mathbb{S}^n$ that contain $A$. For $0 \leq k \leq n$, a $k$-sphere $S$ in $\mathbb{S}^n$ is a \linebreak $k$-dimensional great sub-sphere of $\mathbb{S}^n$. Clearly, every $k$-sphere $S$ is convex.

For $e, v \in \mathbb{S}^n$, the \emph{(geodesic) reflection} of $v$ about $e$ is given by
\[  v^e := -v + 2(v\cdot e)e.\]
A subset $A \subseteq \mathbb{S}^n$ is called centrally-symmetric with center $e\in \mathbb{S}^n$, if $A^e = A$. Let $\mathcal{K}_c(\mathbb{S}^n)$ denote the subspace of $\mathcal{K}(\mathbb{S}^n)$ of all centrally-symmetric spherically convex bodies in $\mathbb{S}^n$. Clearly, if $K \in \mathcal{K}_c(\mathbb{S}^n)$ has center $e \in \mathbb{S}^n$, then $K \subseteq \mathbb{S}_e^+$. Moreover, if $K$ is proper, we have $K \subseteq \mathrm{int}\,\mathbb{S}_e^+$.

The following lemma contains a few useful properties of the spherical Hausdorff metric which we require later on.

\begin{lem}\label{2:lemSpherHD} For $m \in \mathbb{N}$, let $C_m,C \subseteq \mathbb{S}^n$ be closed and $K, L \in \mathcal{K}(\mathbb{S}^n)$ such that
$\delta_s(K, L)<\frac{\pi}{2}$. Then the following statements hold:
\begin{enumerate}
\item[(a)] The sequence $C_m$ converges to $C$ in the spherical Hausdorff metric if and only if it does so in the Hausdorff metric of the ambient space $\mathbb{R}^{n+1}$;
\item[(b)] $\delta_s(K,L) = \delta_s(\mathrm{bd}\, K, \mathrm{bd}\, L)$;
\item[(c)] $\delta_s(K,L) = \delta_s(K^*, L^*)$.
\end{enumerate}
\end{lem}
\noindent {\it Proof.} Statement (a) is a consequence of the fact that $\|u - v\| \leq d_s(u, v) \leq \frac{\pi}{2} \|u - v\|$ for all $u, v \in \mathbb{S}^n$, that is, of the
equivalence of the spherical and the Euclidean distance in the ambient space.

In order to see (b), we use that for $d_s(x, K) < \frac{\pi}{2}$, there exists a unique point $p(K, x)$ in $K$ such that $d_s(x, p(K, x)) < d_s(x, y)$ for all $y\in K$. From this, (b) follows by the same argument as in the linear setting (see, e.g., \textbf{\cite[\textnormal{Lemma 1.8.1}]{Schneider:2014}}).

Finally, a proof of (c) was, for example, given in \textbf{\cite[\textnormal{Hilfssatz 2.2}]{Glasauer:1996}}. \hfill $\blacksquare$

\vspace{0.3cm}

We turn now to one of the most important tools in spherical convexity, the gnomonic projection. First, recall that for $e \in \mathbb{S}^n$, the \emph{gnomonic projection} with respect to $e$ is given by
\[g_e: \mathrm{int}\,\mathbb{S}^+_e \rightarrow \mathbb{R}^n_e, \qquad g_e(v)= \frac{v}{e \cdot v}-e.  \]

In the following, we write $\varrho_e: \mathbb{R}^{n+1} \rightarrow \mathbb{R}^{n+1}$ for the orthogonal reflection about $\mathbb{R}^n_e$ in $\mathbb{R}^{n+1}$ and we let $\mathcal{K}(\mathrm{int}\,\mathbb{S}_e^+)$ and $\mathcal{K}(\mathbb{R}^n_e)$ denote the spaces of convex bodies in $\mathrm{int}\,\mathbb{S}_e^+$ and $\mathbb{R}^n_e$ endowed with the respective Hausdorff metrics. Moreover, recall that for $K \in \mathcal{K}(\mathbb{R}^n_e)$ containing the origin in its interior, the \emph{polar body} of $K$ is defined by $K^\circ = \{x \in \mathbb{R}^n_e: x \cdot y \leq 1 \mbox{ for all } y \in K\}$.
Our next lemma contains several well known properties of the gnomonic projection, the proofs of which can be found, e.g., in \textbf{\cite{Besau:Schuster:2016}} and \textbf{\cite{Gao:Hug:Schneider:2003}}.

\begin{lem}\label{2:propGnom}
For $e \in \mathbb{S}^n$, the gnomonic projection $g_e: \mathrm{int}\,\mathbb{S}_e^+ \to \mathbb{R}^n_e$ has the following properties:
\begin{enumerate}
\item[(a)] The map $g_e$ is a bijection with inverse given by
\[g_e^{-1}: \mathbb{R}^n_e \to \mathrm{int}\,\mathbb{S}^+_e, \qquad g_e^{-1}(x)= \frac{x + e}{\|x + e\|}.\]

\item[(b)] If $S \subseteq \mathbb{S}^n$ is a $k$-sphere, $0 \leq k \leq n - 1$, such that $S \cap \mathrm{int}\,\mathbb{S}_e^+$ is non-empty, then $g_e(S \cap \mathrm{int}\,\mathbb{S}_e^+)$ is a
$k$-dimensional affine subspace of $\mathbb{R}^n_e$. Conversely, $g_e^{-1}$ maps $k$-dimensional affine subspaces of $\mathbb{R}^n_e$ to $k$-hemispheres in $\mathrm{int}\,\mathbb{S}_e^+$.

\item[(c)] The map $g_e$ induces a homeomorphism between $\mathcal{K}(\mathrm{int}\,\mathbb{S}_e^+)$ and $\mathcal{K}(\mathbb{R}^n_e)$.

\item[(d)] For every $u \in \mathrm{int}\,\mathbb{S}^+_e$, we have $g_e(u^e) = -g_e(u)$.

\item[(e)] For every $K \in \mathcal{K}(\mathrm{int}\,\mathbb{S}_e^+)$ containing $e$ in its interior, $g_e(\varrho_e K^*) = g_e(K)^\circ$.
\end{enumerate}
\end{lem}

\pagebreak

For our purposes it is important to know the push-forwards of certain measures on $\mathbb{S}^n$ under gnomonic projection. These are the content of the following lemma.

\begin{lem}\label{2:pushMeasGnom}
Let $e\in\S^n$ and let $g_e\colon\INT\S_e^+\to\R^n_e$ be the gnomonic projection.
\begin{enumerate}
\item[(a)] The push-forward $\widehat{\sigma} := g_e\# \sigma$ under $g_e$ of spherical Lebesgue measure $\sigma$ is absolutely continuous with
density given by $\xi(x) = (1 + \|x\|^2)^{-\frac{n+1}{2}}$.

\item[(b)] For $u\in\mathrm{int}\,\mathrm{S}_e^+$, we have $u\cdot e = \phi(g_e(u))$, where $\phi(x) = (1 + \|x\|^2)^{-\frac{1}{2}}$.

\item[(c)] The push-forward $\widehat{\tau} := g_e\#\tau$ under $g_e$ of the absolutely continuous measure $\tau$ on $\S_e^+$ with density $u \mapsto e \cdot u$ is also absolutely continuous with
density given by $\psi(x) = (1 + \|x\|^2)^{-\frac{n+2}{2}}$.
\end{enumerate}
\end{lem}

\noindent {\it Proof.} In order to see (a), we need the Jacobian of the inverse $g_e^{-1}$ at $x \in \mathbb{R}^n_e$. It was, for example, calculated in
\textbf{\cite[\textnormal{Proposition 4.2}]{Besau:Werner:2016}} and is given by
\[Jg_e^{-1}(x) = (1 + \|x\|^2)^{-\frac{n+1}{2}}=\xi(x).\]
Thus, by the area formula (see e.g., \textbf{\cite[\textnormal{Theorem 8.1}]{Maggi:2012}}), we have
\[ \sigma(A) = \int_{g_e(A)} (1 + \|x\|^2)^{-\frac{n+1}{2}}\, dx,\]
for every Borel set $A \subset \mathrm{int}\,\mathbb{S}_e^+$, which proves (a).

Statement (b) follows from Pythagoras' theorem, since $\|g_e(u)\|^2 + 1 = (u\cdot e)^{-2}$, and, finally, combining (a) and (b) yields statement (c). \hfill $\blacksquare$

\vspace{0.3cm}

Next, we discuss the notion of centroids for certain subsets of the unit sphere which we use in this paper (for other notions of centroids on $\mathbb{S}^n$, cf.\ \textbf{\cite{Galperin:1993}}).

\begin{defi}
For $\{u_1,\ldots, u_N\} \subseteq \mathbb{S}^n$ and a Borel subset $A \subseteq \mathbb{S}^n$ such that $\sigma(A)>0$, we define their respective \emph{spherical centroids} by
\[  c_s(u_1,\ldots, u_N) := \frac{\sum_{i=1}^N u_i}{\left\|\sum_{i=1}^N u_i\right\|} \qquad \textrm{and} \qquad c_s(A):= \frac{\int_A u\, d\sigma(u)}{\left\|\int_A u\, d\sigma(u)\right\|}\]
whenever they exist, that is, whenever the denominators are non-zero.
\end{defi}

While this definition of spherical centroids makes use of the vector space structure of the ambient space, it is well known that both $c_s(u_1,\ldots, u_N)$ and $c_s(A)$ can also be defined (with more complicated formulae) intrinsically, that is, by making use only of the metric structure of the sphere (see e.g., Galperin \textbf{\cite{Galperin:1993}}, where he also characterized $c_s$ by a set of natural properties).

In order to carry out explicit computations later on, we combine now the gnomonic projection with spherical centroids. To this end, we need to consider centroids in $\mathbb{R}^n$ with respect to arbitrary densities.

\begin{defi}
For $\{x_1, \ldots, x_N\} \subseteq \mathbb{R}^n$ and a positive function $f: \mathbb{R}^n \to \mathbb{R}^+$, let
\begin{equation} \label{defcf}
c_f(x_1,\dots,x_N) := \frac{1}{\sum_{i=1}^N f(x_i)}\sum_{i=1}^N f(x_i)x_i.
\end{equation}
For an absolutely continuous measure $\mu$ on $\mathbb{R}^n$ and a bounded Borel subset $A \subseteq \mathbb{R}^n$ such that $\mu(A)>0$, we define the \emph{$\mu$-centroid} of $A$ by
\begin{equation} \label{defcmu}
c_\mu(A) := \frac{1}{\mu(A)} \int_A x\, d\mu(x).
\end{equation}
\end{defi}

Our next lemma is critical for the proof of relation (\ref{gnomGammas}) and Theorem \ref{1:thmApprox}. Here and in the following, we use again the notation from Lemma \ref{2:pushMeasGnom}.

\begin{lem}\label{2:csGnom}
If $\{u_1,\ldots,u_N\} \subseteq \mathrm{int}\,\mathbb{S}_e^+$ for some $e \in \mathbb{S}^n$, then
\begin{equation} \label{gnomcent1}
g_e(c_s(u_1,\ldots,u_N)) = c_\phi(g_e(u_1),\dots,g_e(u_N)).
\end{equation}
If $A \subseteq \mathrm{int}\,\mathbb{S}_e^+$ is a Borel subset such that $\sigma(A)>0$, then
\begin{equation} \label{gnomcent2}
g_e(c_s(A)) = c_{\,\widehat{\tau}}(g_e(A)).
\end{equation}
\end{lem}

\noindent {\it Proof.} By Lemma \ref{2:pushMeasGnom} (b) and definition (\ref{defcf}), relation (\ref{gnomcent1}) follows from
\[   g_e\left(\sum_{i=1}^N u_i\right) = \frac{\sum_{i=1}^N u_i}{\sum_{i=1}^N u_i\cdot e} - e
		 = \frac{\sum_{i=1}^N (u_i\cdot e)\left(\frac{u_i}{u_i\cdot e}-e\right)}{\sum_{i=1}^N u_i\cdot e}
	   = \frac{1}{\sum_{i=1}^N \phi(x_i)}\sum_{i=1}^N \phi(x_i)x_i,\]
where $x_i=g_e(u_i)$, $1 \leq i \leq N$.
In order to prove (\ref{gnomcent1}), we use again the area formula (see e.g., \textbf{\cite[\textnormal{Remark 8.3}]{Maggi:2012}}) to obtain
\[  g_e(c_s(A)) = \frac{\int_A u\, d\sigma(u)}{e \cdot \int_A u\, d\sigma(u)} - e
	   = \frac{\int_{g_e(A)} g_e^{-1}(x)Jg_e^{-1}(x)\, dx}{\int_{g_e(A)} e \cdot g_e^{-1}(x)Jg_e^{-1}(x)\, dx} - e.\]
Since for $x \in \mathbb{R}^n_e$, we have $\|x + e\|^2 = 1 + \|x\|^2$, and by the proof of Lemma \ref{2:pushMeasGnom} (a) $Jg_e^{-1}(x) = \xi(x)$, we conclude that
\begin{align*}
  g_e(c_s(A)) = \frac{\int_{g_e(A)} (x + e) \xi(x)\, dx}{\int_{g_e(A)} e \cdot (x + e) \xi(x)\, dx} - e
	&=  \frac{\int_{g_e(A)} (x + e) \, d\phantom{.}\!\widehat{\tau}(x)}
	     {\int_{g_e(A)} 1 \, d\phantom{.}\!\widehat{\tau}(x)} - e \\
	&=  \frac{1}{\widehat{\tau}(g_e(A))} \int_{g_e(A)} x \, d\phantom{.}\!\widehat{\tau}(x)
	 =  c_{\,\widehat{\tau}}(g_e(A)).
\end{align*}

\vspace{-0.6cm}

\hfill $\blacksquare$

\vspace{0.3cm}

We conclude this section by collecting a number of properties of spherical centroids for quick later reference.

\begin{prop}\label{2:propSpherCent}
Let $\{u_1,\ldots, u_N\} \subseteq \mathbb{S}^n$ and $K_m, K\in\K(\S^n)$, $m\in\mathbb{N}$, such that their spherical centroids exist. Then the map $c_s$ has the following properties:
\begin{enumerate}
\item[(a)] It is continuous, that is, if $u_{i,m}\to u_i$ for $1\leq i\leq N$ and $K_m\to K$, $m\in\N$, in the spherical Hausdorff metric, then
\[ c_s(u_{1,m},\ldots, u_{N,m})\to c_s(u_1,\ldots, u_N) \qquad \textrm{and} \qquad  c_s(K_m)\to c_s(K).\]

\item[(b)] It is $\ORTHO(n+1)$-equivariant, that is, for every $\vartheta \in \mathrm{O}(n+1)$, we have
\[ c_s(\vartheta u_1,\ldots, \vartheta u_N) = \vartheta c_s(u_1,\ldots, u_N) \qquad \textrm{and} \qquad c_s(\vartheta K) = \vartheta c_s(K).\]
				
\item[(c)] It is proper, that is, $c_s(K)\in\mathrm{int}\, K$.

\item[(d)] It is consistent, that is, if $U_1,\ldots, U_N$ are independent random variables uniformly distributed in $K$, then
\[ c_s(U_1,\ldots, U_N) \to c_s(K)\]
almost surely as $N$ tends to infinity.
\end{enumerate}
\end{prop}
\noindent {\it Proof.} Property (a) is trivial in the discrete case and follows for convex bodies from the continuity of spherical volume in the Hausdorff topology on convex bodies (see e.g., \textbf{\cite[\textnormal{Hilfssatz 2.4}]{Glasauer:1996}}) since
\[  \left|\int_{K_m} u \, d\sigma(u) - \int_K u \, d\sigma(u)\right|  \leq \sigma(K_m \triangle K)\to 0.\]

Property (b) is also trivial in the discrete case and for convex bodies a simple consequence of the $\mathrm{O}(n+1)$-invariance of spherical Lebesgue measure and the transformation rule for integrals.

In order to see (c), note that $u \in \mathrm{int}\, K$ if and only if $w\cdot u<0$ for all $w\in K^*$. Now since $w\cdot c_s(K)<0$ for all $w\in K^*$, by definition, we obtain the desired property.

Finally, since $\sigma(\mathrm{bd}\, K)=0$ and $\mathrm{int}\, K$ is proper, we may assume for the proof of (d) that $K \subset \mathrm{int}\,\mathbb{S}_e^+$ for some $e\in\S^n$. Then, by Lemma \ref{2:pushMeasGnom} (a), the random variables $X_i := g_e(U_i)$, $1\leq i\leq N$, are independent and identically distributed according to
\[ \frac{\mathbbm{1}_{g_e(K)}}{\widehat{\sigma}(g_e(K))}\widehat{\sigma}. \]
Moreover, by Lemma \ref{2:csGnom}, we have
\[g_e(c_s(U_1,\ldots, U_N)) = \frac{1}{\sum_{i=1}^N \phi(X_i)} \sum_{i=1}^N \phi(X_i) X_i.\]
But, by the strong law of large numbers (see e.g., \textbf{\cite[\textnormal{Theorem 8.3.5}]{dudley2002}}),
\[\sum_{i=1}^N \phi(X_i) \to \int_{g_e(K)} \phi(x)\, d\phantom{.}\!\widehat{\sigma}(x) = \widehat{\tau}(g_e(K))\]
and
\[\sum_{i=1}^N \phi(X_i)X_i \to \int_{g_e(K)} \phi(x)x_i \, d\phantom{.}\!\widehat{\sigma}(x) = \int_{g_e(K)} x\, d\phantom{.}\!\widehat{\tau}(x)\]
almost surely as $N\to\infty$. Since, by the continuous mapping theorem, the product of almost surely convergent sequences of random variables converges almost surely to the product of their limits, we obtain from another application of Lemma \ref{2:csGnom},
\[g_e(c_s(U_1,\ldots, U_N)) \to \frac{1}{\widehat{\tau}(g_e(K))} \int_{g_e(K)}x\, d\phantom{.}\!\widehat{\tau}(x) = c_{\,\widehat{\tau}}(g_e(K)) = g_e(c_s(K)),\]
almost surely as $N\to\infty$, which, by Lemma \ref{2:propGnom} (c), yields property (d). \hfill $\blacksquare$

\vspace{1cm}

\centerline{\large{\bf{ \setcounter{abschnitt}{3}
\arabic{abschnitt}. Centroid bodies}}}

\reseteqn \alpheqn \setcounter{theorem}{0}

\vspace{0.6cm}

In the first part of this section we discuss the definition and properties of weighted centroid bodies in linear vector spaces. The second part is devoted to spherical centroid bodies and their basic properties. We also establish a few auxiliary results required for the proof of Theorem \ref{1:thmApprox}.

We begin with a definition of weighted centroid bodies of arbitrary convex bodies in $\mathbb{R}^n$, different from the one for origin-symmetric bodies given in the introduction (Lemma \ref{2:muCentBound} below shows that the two definitions coincide on origin-symmetric convex bodies). To this end, recall that a convex body $K \in \mathcal{K}(\mathbb{R}^n)$ is uniquely determined by the values of its support function $h(K,u)=\max\{u \cdot x: x \in K\}$, $u \in \mathbb{S}^{n-1}$, and that every even, positive, and sublinear function on $\mathbb{R}^n$ is the support function of an origin-symmetric convex body in $\mathbb{R}^n$ (see e.g., \textbf{\cite[\textnormal{Theorem 1.7.1}]{Schneider:2014}}).

\begin{defi}
For $\{x_1,\ldots,x_N\} \subseteq \mathbb{R}^n$ and a positive function $f: \mathbb{R}^n \to \mathbb{R}^+$, define
\begin{equation} \label{fcentx1xN}
h(\Gamma_f(x_1,\ldots,x_N), u) := \frac{1}{\sum_{i=1}^N f(x_i)} \sum_{i=1}^N |u \cdot f(x_i)x_i|.
\end{equation}
For a finite Borel measure $\mu$ on $\mathbb{R}^n$ with positive density and $L \in \mathcal{K}(\mathbb{R}^n)$, define the \emph{$\mu$-centroid body} of $L$ by
\begin{equation} \label{defmucent17}
h(\Gamma_\mu L, u) := \frac{1}{\mu(L)} \int_L |u \cdot y|\, d\mu(y).
\end{equation}
\end{defi}

Note that, by our assumption on $\mu$, $\Gamma_\mu L$ is an origin-symmetric convex body for every $L \in \mathcal{K}(\mathbb{R}^n)$. While $\Gamma_f(x_1,\ldots, x_N)$ is, in general, always an origin-symmetric, compact, convex set, it has non-empty interior if and only if $\mathrm{span}\,\{x_1,\ldots,x_N\}=\mathbb{R}^n$. It is also worth noting that when $\mu$ is taken to be Lebesgue measure, (\ref{defmucent17}) defines Blaschke's classical centroid body (of the not necessarily origin-symmetric) body $L$. In the following, when $f \equiv 1$ in (\ref{fcentx1xN}), we simply write
$\Gamma(x_1,\ldots,x_N)$ and use $h([-z,z],u)=|u\cdot z|$ for every $z \in \mathbb{R}^n$, to see that, in this case,
\[  \Gamma(x_1,\ldots,x_N) = \frac{1}{N}\sum_{i=1}^N [-x_i, x_i].\]

Our first goal is to relate weighted centroid bodies with weighted centroids. In the discrete case this is the content of the following easy lemma.

\begin{lem}\label{2:fCentroidBodies}
Let $\{x_1,\ldots,x_N\} \subseteq \mathbb{R}^n$ be a finite subset and assume that $f: \mathbb{R}^n \to \mathbb{R}^+$ is even. Then
\[ \Gamma_f(x_1,\ldots, x_N) = \mathrm{conv}\left\{ c_f(\pm x_1, \ldots, \pm x_N) \right\}.\]
\end{lem}
\noindent {\it Proof.} Since for arbitrary $z_1,\ldots, z_N \in \mathbb{R}^n$, we have
\[\sum_{i=1}^N [-z_i, z_i] = \mathrm{conv}\{\pm z_1 \pm \cdots \pm z_N \},\]
we obtain
\[\Gamma_f(x_1,\ldots, x_N) = \mathrm{conv}\left\{\frac{\pm f(x_1)x_1 \pm \cdots \pm f(x_N)x_N}{\sum_{i=1}^N f(x_i)} \right\}.\]
But, since $f$ is even, this is equal to $\mathrm{conv}\left\{ c_f(\pm x_1, \ldots, \pm x_N) \right\}$. \hfill $\blacksquare$

\vspace{0.3cm}

In contrast to $\Gamma_f(x_1, \ldots, x_N)$, which is, as a Minkowski sum of line segments, always a polytope, our next lemma shows that the boundary of $\Gamma_\mu L$ exhibits higher regularity.
For the classical centroid body this was first proved by Petty \textbf{\cite{Petty:1961}} and we follow his arguments closely (see also, \textbf{\cite[\textnormal{Theorem 1.2}]{Huang:Slomka:Werner:2018}} for a recent variant). In order to state the result precisely, recall that a convex body $L$ is said to be of class $C^2_+$ if the boundary of $L$ is a $C^2$ submanifold of $\mathbb{R}^n$ with everywhere positive Gau\ss\ curvature.

\begin{lem}\label{2:lemMuCentBodC2}
Let $\mu$ be a finite Borel measure on $\mathbb{R}^n$ with positive bounded density and $L \in \mathcal{K}(\mathbb{R}^n)$. Then $\Gamma_\mu L$ is of class $C^2_+$. In particular, it is strictly convex.
\end{lem}
\noindent {\it Proof.} We first want to show that $h(\Gamma_\mu L, \cdot)$ is twice differentiable. To this end, we compute its directional derivative at $x \in \mathbb{R}^n$ in the direction
$u \in \mathbb{S}^{n-1}$ by
\[\lim_{t\to 0^+} \frac{h(\Gamma_\mu L, x + tu) - h(\Gamma_\mu L, x)}{t} = \frac{1}{\mu(L)} \left(\int_{L\cap H_x^+}\!\! u \cdot y\, d\mu(y) - \int_{L\cap H_x^-}\!\! u \cdot y\, d\mu(y) \right).\]
Consequently, the gradient of $h(\Gamma_\mu L, \cdot)$ exists and is given by
\[\nabla h(\Gamma_\mu L, \cdot)(x) = \frac{1}{\mu(L)} \left(\int_{L\cap H_x^+} y\, d\mu(y) - \int_{L\cap H_x^-} y\, d\mu(y) \right). \]
In order to compute second derivatives at $\bar{x} \in \mathbb{R}^n$, we choose an orthonormal coordinate frame $\{e_1,\ldots, e_n\}$ such that $\bar{x}=(0,\ldots,0, \bar{x}_n)^{\mathrm{T}}$, where $\bar{x}_n>0$ (see e.g., \textbf{\cite[\textnormal{p.\ 57}]{Busemann:1958}}). Since $\nabla h(\Gamma_\mu L, \cdot)(x)$ is $0$-homogeneous in $x$, it follows that
\[\frac{\partial^2 h(\Gamma_\mu L, \cdot)}{\partial e_i\partial e_n}(\bar{x}) = 0\]
for $1\leq i\leq n$. Letting $x=(0,\ldots, 0, x_j, 0,\ldots, 0, \bar{x}_n)^{\mathrm{T}}$ for $j<n$, we get for $i,j < n$,
\begin{equation} \label{smoothproof1}
\frac{\frac{\partial h(\Gamma_\mu L, \cdot)}{\partial e_i}(x) - \frac{\partial h(\Gamma_\mu L, \cdot)}{\partial e_i}(\bar{x})}{x_j}
				= \frac{2}{x_j \mu(L)}\! \left(\int_{L\cap H_x^+ \cap H_{\bar{x}}^-}\!\!\!\! y_i\, d\mu(y) -\! \int_{L\cap H_x^- \cap H_{\bar{x}}^+}\!\!\!\! y_i\, d\mu(y) \right)\!.
\end{equation}
In order to compute the limit $x_j \to 0$ in (\ref{smoothproof1}), we make the change of variables $y_1 = v_1$, $\ldots$, $y_{n-1} = v_{n-1}$, and $y_n = v_j \tan v_n$.
The Jacobian of this transformation is given by $J(v) = v_j\sec^2 v_n$. Note that it is negative on $L\cap H_x^- \cap H_{\bar{x}}^+$ when $x_j>0$ and on $L\cap H_x^+ \cap H_{\bar{x}}^-$ when $x_j<0$.
Letting $\alpha := \arctan(|x_j|/\bar{x}_n)$,
\[H(v_n):= \left\{ \begin{array}{lc} y_n = (\tan v_n)y_j & \mbox{ for } x_j<0, \\ y_n = -(\tan v_n)y_j & \mbox{ for } x_j>0, \end{array} \right.\]
and $L^\pm(v_n):= L \cap H(v_n) \cap H_{e_j}^\pm$ for $v_n>0$, the right hand side of (\ref{smoothproof1}) becomes
\begin{equation} \label{smoothproof2}
  \frac{2}{\bar{x}_n\mu(L)\tan\alpha} \left( \int_0^\alpha \sec^2 v_n\, \varphi^+(v_n)\, dv_n + \int_0^\alpha\sec^2 v_n \, \varphi^-(v_n)\, dv_n \right)\!,
\end{equation}
where
\[\varphi^\pm(s) := \int_{L^\pm(s)} v_i v_j\, f_\mu(v_1,\ldots,v_{n-1},v_j \tan s)\,  dv_1\cdots dv_{n-1}\]
with $f_\mu$ being the density of $\mu$. In order to compute the limit $\alpha\to 0$ in (\ref{smoothproof2}), we use that, by the mean value theorem, for every function $\zeta$ which is continuously differentiable near $0$ such that $\zeta(0) = 0$ and every $\varphi$ continuous near $0$,
\[\lim_{\alpha \to 0}\frac{1}{\zeta(\alpha)}\int_0^\alpha \zeta'(s) \varphi(s)\, ds = \varphi(0).\]
Taking here $\zeta(s) := \tan s$, $\varphi = \varphi^\pm$, and letting $\alpha \to 0$ in (\ref{smoothproof2}) as well as changing back to the variables $y_1, \ldots, y_{n-1}$, we obtain
\[  \lim_{x_j\to 0} \frac{\frac{\partial h(\Gamma_\mu L, \cdot)}{\partial e_i}(x) - \frac{\partial h(\Gamma_\mu L, \cdot)}{\partial e_i}(\bar{x})}{x_j}
	= \frac{2}{\bar{x}_n \mu(L)}\! \int_{L\cap H_{\bar{x}}}\!\! y_i y_j\, f_\mu(y_1,\ldots, y_{n-1}, 0)\, dy_1\cdots dy_{n-1},\]
provided we can show that $\varphi^\pm$ is continuous near 0. But since $L\subseteq B_R$ for some Euclidean ball $B_R$ of radius $R$ in $\mathbb{R}^n$ and for every $s_0 \in (0,\varepsilon)$, $L^\pm(s)\to L^\pm(s_0)$ in the Hausdorff metric in $\mathbb{R}^n$ as $s\to s_0$, it is not difficult to see that
\[|\varphi^{\pm}(s)-\varphi^{\pm}(s_0)| \to 0. \]

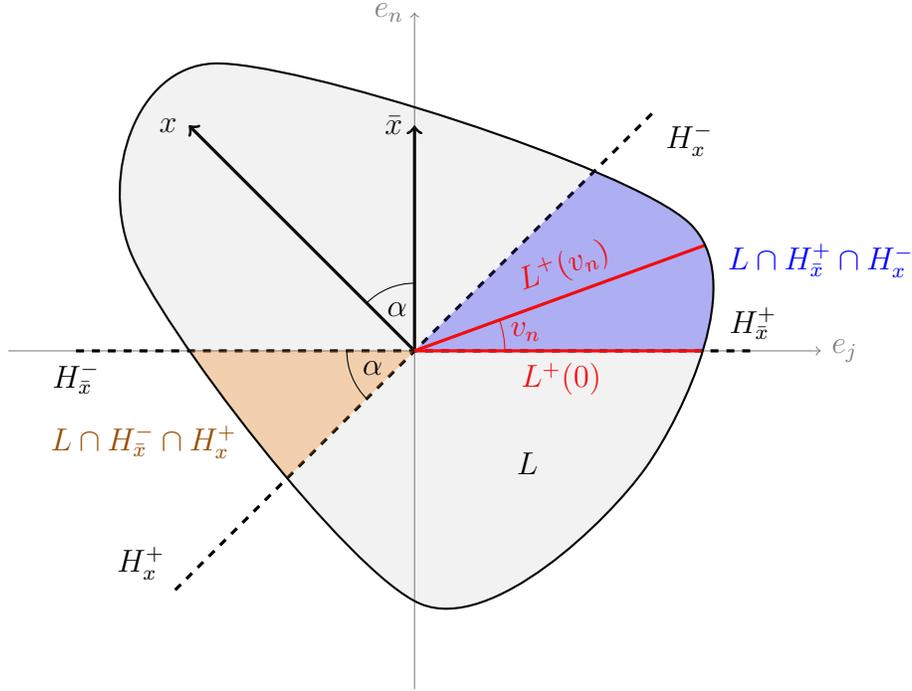
\begin{figure}[h]
\centering
\begin{tikzpicture}[scale=3]
\draw[->, gray] (-1.8, 0) -- (1.8, 0) node[right] {$e_j$};
\draw[->, gray] ( 0,-1.5) -- ( 0, 1.5) node[left] {$e_n$};
\draw[->, very thick] (0, 0) -- ( 0, 1) node[left] {$\bar{x}$};
\draw[->, very thick] (0, 0) -- (-1, 1) node[left] {$x$};
\draw[very thick, dashed] (-1.5, 0) -- (1.5, 0)
node[above] {$H_{\bar{x}}^+$};
\node[below] at (-1.5, 0) {$H_{\bar{x}}^-$};
\draw[very thick, dashed] ({-1.5 / sqrt(2)}, {-1.5 / sqrt(2)}) -- ({1.5 / sqrt(2)}, {1.5 / sqrt(2)})
node[below right] {$H_x^-$};
\node[above left] at ({-1.5 / sqrt(2)}, {-1.5 / sqrt(2)}) {$H_x^+$};
\begin{scope}
\clip[rotate=170, scale=1.1] plot[smooth cycle]
coordinates {(-1, 0.3) (-1, -0.7) (1, -1) (1.2, -0.2) (-0.2, 1)};
\fill[blue, opacity=0.3] ( 1.5, 0) arc( 0: 45:1.5) -- (0, 0) -- cycle;
\fill[orange, opacity=0.3] (-1.5, 0) arc(180:225:1.5) -- (0, 0) -- cycle;
\draw[red, very thick] (0, 0) -- (20:1.5) node[midway, above, rotate=20] {$L^+(v_n)$};
\draw[red] (0.4, 0) arc(0:20:0.4);
\node[red] at (10:0.5) {$v_n$};
\draw[red, very thick] (0, 0) -- (0:1.3) node[midway, below] {$L^+(0)$};
\end{scope}
\draw[thick, rotate=170, scale=1.1] plot[smooth cycle]
coordinates {(-1, 0.3) (-1, -0.7) (1, -1) (1.2, -0.2) (-0.2, 1)};
\fill[gray, rotate=170, scale=1.1, opacity = 0.1] plot[smooth cycle]
coordinates {(-1, 0.3) (-1, -0.7) (1, -1) (1.2, -0.2) (-0.2, 1)};
\node[blue] at (1.8, 0.4) {$L\cap H_{\bar{x}}^+ \cap H_x^-$};
\node[orange!60!black] at (-1.2, -0.4) {$L\cap H_{\bar{x}}^- \cap H_x^+$};
\draw ( 0, 0.3) arc( 90:135:0.3);
\node at (112.5:0.2) {$\alpha$};
\draw (-0.3, 0) arc(180:225:0.3);
\node at (202.5:0.2) {$\alpha$};
\node at (0.5, -0.5) {$L$};
\end{tikzpicture}
\caption{Sketch for the proof of Lemma 3.2}
\label{fig:1}
\end{figure}

Letting $A:=(h_{ij})_{i,j=1}^{n-1}$ denote the Hessian matrix of $h(\Gamma_\mu L, \cdot)$ at $\bar{x} \in \mathbb{R}^n$ (w.r.t.\ $\{e_1,\ldots,e_{n-1}\}$), we can now conclude that
for any $b \in \mathbb{R}^{n-1}\setminus\{0\}$,
\[ b \cdot Ab = \frac{2}{\bar{x}_n \mu(L)} \int_{L\cap H_{\bar{x}}} (b\cdot y)^2\,  f_\mu(y_1,\ldots, y_{n-1}, 0)\, dy_1\cdots dy_{n-1} > 0,\]
that is, $A$ is a positive-definite matrix. Since $\bar{x}$ was arbitrary, it is well known (cf.\ \textbf{\cite[\textnormal{Section 2.5}]{Schneider:2014}}), that this implies that $\Gamma_\mu L$ is of class $C_+^2$ \hfill $\blacksquare$

\vspace{0.3cm}

The following lemma shows that (\ref{introdefmycent}) and (\ref{defmucent17}) define the same convex bodies in the case of an even measure and an origin-symmetric convex body $L$.

\begin{lem}\label{2:muCentBound}
Let $\mu$ be a finite even Borel measure on $\mathbb{R}^n$ with positive density and assume that $L \in \mathcal{K}(\mathbb{R}^n)$ is origin-symmetric.
Then
\[\mathrm{bd}\, \Gamma_\mu L = \{ c_\mu(L\cap H_u^+): u\in\mathbb{S}^{n-1}\}.\]
\end{lem}
\noindent {\it Proof.} Since $\mu$ is even and $L$ origin-symmetric, we have for $u \in \mathbb{S}^{n-1}$,
\begin{align*}
  h(\Gamma_\mu L, u) &= \frac{1}{\mu(L)}\int_L |u\cdot y|\, d\mu(y) = \frac{2}{\mu(L)}\int_{L\cap H_u^+} u\cdot y\, d\mu(y) \\
	   &= u \cdot \frac{1}{\mu(L\cap H_u^+)}\int_{L\cap H_u^+} y\, d\mu(y)
		 = u\cdot c_\mu(L\cap H_u^+).
\end{align*}
Thus, by the definition of support functions, $c_\mu(L\cap H_u^+)\in \mathrm{bd}\,\Gamma_\mu L$. Since $\Gamma_\mu L$ is strictly convex by Lemma \ref{2:lemMuCentBodC2}, all boundary points are obtained in this way. \hfill $\blacksquare$

\vspace{0.3cm}

We also note that if $L$ in Lemma \ref{2:muCentBound} is not origin-symmetric, then a similar computation shows that every boundary point of $\Gamma_\mu L$ is a convex combination of $c_\mu(L\cap H_u^+)$ and $-c_\mu(L\cap H_u^-)$ for $u \in \mathbb{S}^{n-1}$ (cf.\ \textbf{\cite[\textnormal{Section 9.1}]{Gardner:2006}}).

\vspace{0.3cm}

We now turn our focus towards spherical centroid bodies and first recall their definition:
For $K \in \mathcal{K}_c(\mathbb{S}^n)$ with center $e \in \mathbb{S}^n$, its spherical centroid body is defined by
\[\mathrm{bd}\,\Gamma_s K := \{c_s(K\cap\mathbb{S}_u^+):u\in\mathbb{S}_e\}.\]
Note that since $K$ has nonempty interior, $c_s(K\cap\mathbb{S}_u^+)$ exists for every $u \in \mathbb{S}_e$.
Moreover, since $K \subseteq \mathbb{S}_e^+$, Proposition \ref{2:propSpherCent} (c) implies that $\mathrm{bd}\,\Gamma_s K$ is contained in $\mathrm{int}\,\mathbb{S}_e^+$ and, hence, we can consider its gnomonic projection.

\begin{prop}\label{3:scbGnom}
Let $K \in \mathcal{K}_c(\mathbb{S}^n)$ have center $e \in \mathbb{S}^n$ and let $g_e: \mathrm{int}\,\mathbb{S}_e^+ \to \mathbb{R}^n_e$ denote the gnomonic projection.
Then
\[ g_e(\mathrm{bd}\,\Gamma_s K) = \mathrm{bd}\,\Gamma_{\widehat{\tau}} \,g_e(K \cap \mathrm{int}\,\mathbb{S}_e^+).\]
\end{prop}
\noindent {\it Proof.} Let us first assume that $K$ is proper, that is, $K = K \cap \mathrm{int}\,\mathbb{S}_e^+$. Then, by Lemma \ref{2:csGnom}, we have
\[g_e(\mathrm{bd}\,\Gamma_s K) = \{c_{\widehat{\tau}}(g_e(K\cap \S_u^+)): u \in \mathbb{S}_e\}. \]
But since $g_e(K \cap \mathbb{S}_u^+) = g_e(K) \cap H_u^+$ for every $u \in \mathbb{S}_e$, $\widehat{\tau}$ is even, and $g_e(K)$ is origin-symmetric, it follows from Lemma \ref{2:muCentBound} that
\[g_e(\mathrm{bd}\,\Gamma_s K) = \{c_{\widehat{\tau}}(g_e(K) \cap H_u^+): u\in \mathbb{S}_e\} = \mathrm{bd}\,\Gamma_{\widehat{\tau}}\, g_e(K).\]

Now, if $K$ is not proper, then $g_e(K \cap \mathrm{int}\,\mathbb{S}_e^+)$ is a closed, convex, and origin-symmetric set in $\mathbb{R}^n_e$ with nonempty interior which is unbounded.
However, since $\widehat{\tau}$ has finite first moments, that is,  $\int_{\mathbb{R}^n} y_i\, d\phantom{.}\!\widehat{\tau}(y)<\infty$ for all $1\leq i\leq n$, (\ref{defmucent17}) still makes sense and defines a convex body $\Gamma_{\widehat{\tau}}\, g_e(K \cap \mathrm{int}\,\mathbb{S}_e^+)$ in $\mathbb{R}^n_e$.

\pagebreak

Moreover, since the density function $\psi$ of $\widehat{\tau}$ is radially symmetric, radially decreasing (see Section 4 for definitions), and satisfies
\[\int_{\mathbb{R}^{n-1}} y_i y_j\, \psi(y_1,\ldots, y_{n-1}, 0)\, dy_1\cdots dy_{n-1}<\infty\]
for all $1\leq i,j \leq n-1$, it is not difficult to show that Lemma \ref{2:lemMuCentBodC2} also holds for $\Gamma_{\widehat{\tau}}\, g_e(K \cap \mathrm{int}\,\mathbb{S}_e^+)$ (the key is to prove continuity near 0 of the respective functions $\varphi^{\pm}$ from (\ref{smoothproof2})), and, therefore, so does Lemma \ref{2:muCentBound}. Consequently, since removing sets of measure zero does not affect centroid computations, the arguments from the first part of the proof yield the desired relation also for nonproper $K$. \hfill $\blacksquare$

\vspace{0.3cm}

As a consequence of Proposition \ref{3:scbGnom} and the properties of the gnomonic projection, we obtain that the spherical centroid body map is, in fact, well defined.

\begin{koro}\label{3:scbGnom2}
Let $K \in \mathcal{K}_c(\mathbb{S}^n)$ have center $e \in \mathbb{S}^n$ and let $g_e: \mathrm{int}\,\mathbb{S}_e^+ \to \mathbb{R}^n_e$ denote the gnomonic projection.
Then $\Gamma_s K \in \mathcal{K}_c(\mathbb{S}^n)$ is proper and has center $e$. Moreover,
\[ g_e(\Gamma_s K) = \Gamma_{\widehat{\tau}\,} g_e(K \cap \mathrm{int}\,\mathbb{S}_e^+).\]
\end{koro}
\noindent {\it Proof.} Since $\Gamma_{\widehat{\tau}} L$ is an origin-symmetric convex body in $\mathbb{R}^n$ for any (possibly unbounded) closed, convex set $L\subseteq \mathbb{R}^n$ with nonempty interior, the statement follows from Proposition \ref{3:scbGnom} and Lemma \ref{2:propGnom}. \hfill $\blacksquare$

\vspace{0.3cm}

With the following proposition we collect several basic properties of the spherical centroid body map.

\begin{prop}\label{3:scbProp}
The spherical centroid body map $\Gamma_s: \mathcal{K}_c(\mathbb{S}^n) \rightarrow \mathcal{K}_c(\mathbb{S}^n)$ has the following properties:
\begin{enumerate}
\item[(a)] It is $\mathrm{O}(n+1)$-equivariant, that is, $\Gamma_s (\vartheta K) = \vartheta\Gamma_s K$ for all $\vartheta \in \mathrm{O}(n+1)$;
\item[(b)] It is continuous;	
\item[(c)] It is injective on bodies of equal spherical volume;
\item[(d)] $\Gamma_s K$ is of class $C^2_+$ for every $K \in \mathcal{K}_c(\mathbb{S}^n)$.
\end{enumerate}
\end{prop}
\noindent {\it Proof.} Property (a) is an immediate consequence of Proposition \ref{2:propSpherCent} (b) and the definition of $\Gamma_s K$.

In order to prove (b), let $K_m, K\in \mathcal{K}_c(\mathbb{S}^n)$ such that $K_m\to K$ in the spherical Hausdorff metric. In addition, let us first also assume that all $K_m$ and $K$ have the same center $e \in \mathbb{S}^n$. In this case, it follows that $K_m \cap \mathbb{S}_u^+ \to K \cap \mathbb{S}_u^+$ for all $u\in \mathbb{S}_e$ and thus, by Proposition \ref{2:propSpherCent} (a),
\begin{equation} \label{gammascont1}
c_s(K_m \cap \mathbb{S}_u^+) \to c_s(K \cap \mathbb{S}_u^+).
\end{equation}
Now note the following two consequences of (\ref{gammascont1}):
\begin{enumerate}
\item[(i)] For every $v = c_s(K \cap \mathbb{S}_u^+) \in \mathrm{bd}\,\Gamma_s K$, there exists a sequence $v_m \in \mathrm{bd}\,\Gamma_s K_m$ such that $v_m\to v$,
namely $v_m = c_s(K_m \cap \mathbb{S}_u^+)$.
\end{enumerate}

\pagebreak

\begin{enumerate}
\item[(ii)] For every convergent subsequence $v_{m_l} \to v$ with $v_{m_l} \in \mathrm{bd}\,\Gamma_s K_m$, we have $v \in \mathrm{bd}\,\Gamma_s K$. Indeed, since
$v_{m_l} = c_s(\S_{u_{m_l}}^+ \cap K_{m_l})$ for some $u_{m_l} \in \mathbb{S}_e$ and $\mathbb{S}_e$ is compact, we find a subsequence (which we again call $u_{m_l}$) such that $u_{m_l} \to u \in \mathbb{S}_e$ and, thus,
\[ v = \lim_{l\to\infty} c_s\left(\S_{u_{m_l}}^+ \cap K_{m_l}\right) = c_s(\SP_u^+ \cap K) \in \BD \Gamma_s K.\]
Moreover, the sequence $\BD\Gamma_s K_m\subset\S^n$ is bounded in $\R^{n+1}$.
\end{enumerate}
It is well known (cf.\ \textbf{\cite[\textnormal{p.\ 69}]{Schneider:2014}}) that (i) and (ii) imply $\mathrm{bd}\,\Gamma_s K_m \to \mathrm{bd}\,\Gamma_s K$ in the Hausdorff metric in $\mathbb{R}^{n+1}$.
Consequently, by Lemma \ref{2:lemSpherHD} (a) and (b), $\Gamma_s K_m \to \Gamma_s K$ in the spherical Hausdorff metric.

It remains to settle the case where the bodies $K_m$ have center $e_m \in \mathbb{S}^n$ and $K$ has center $e \in \mathbb{S}^n$. Clearly, the convergence $K_m \to K$ implies that $e_m \to e$ as $m \to \infty$. In the following, we make use (twice) of the fact that if $\vartheta_m, \vartheta\in \mathrm{O}(n+1)$ and $K_m, K\in \mathcal{K}_c(\mathbb{S}^n)$, then $\vartheta_m \to \vartheta$ and $K_m\to K$ imply $\vartheta_m K_m \to \vartheta K$ in the spherical Hausdorff metric. To apply this, note that since $e_m \to e$, there exists a sequence $\vartheta_m \in \mathrm{O}(n+1)$ such that $\vartheta_m e_m = e$ and $\vartheta_m \to \mathrm{Id}$. Hence, $\vartheta_m K_m \to K$ and, since all $\vartheta_m K_m$ have center $e$, property (a) and the first part of the proof of (b) imply
\[\vartheta_m \Gamma_s K_m = \Gamma_s (\vartheta_m K_m) \to \Gamma_s K.\]
Making use a second time of the above fact, now for $\vartheta_m^{-1}\to \mathrm{Id}$, yields $\Gamma_s K_m \to \Gamma_s K$ which completes the proof of (b).

In order to prove the injectivity property (c), let $K, L \in \mathcal{K}_c(\mathbb{S}^n)$ such that $\Gamma_s K = \Gamma_s L$ and assume w.l.o.g.\ that $\tau(K) \leq \tau(L)$.
Since $\Gamma_s K$ and $\Gamma_s L$ have the same center as $K$ and $L$, respectively, it follows that $K$ and $L$ have the same center, say $e\in \mathbb{S}^n$. Moreover, by using polar coordinates, we have
\[h(\Gamma_{\widehat{\tau}\,} g_e(K \cap \mathrm{int}\,\mathbb{S}_e^+),u) = \int_{\mathbb{S}^{n-1}} \left|u\cdot v\right| \frac{1}{\tau(K)} \int_{0}^{\rho_{g_e(K \cap \mathrm{int}\,\mathbb{S}_e^+)}(v)} \!\!\frac{r^{n}}{(1+r^2)^{\frac{n+2}{2}}}\, d{r} \, d{v}, \]
where $\rho_{g_e(K \cap \mathrm{int}\,\mathbb{S}_e^+)}$ denotes the (possibly infinite) radial function of $K \cap \mathrm{int}\,\mathbb{S}_e^+$. \linebreak Hence, by our assumption that $\Gamma_s K = \Gamma_s L$, Corollary \ref{3:scbGnom2}, and the injectivity of the spherical cosine transform on even functions
(cf.\ \textbf{\cite[\textnormal{Theorem~C.2.1}]{Gardner:2006}}), we conclude that
\[\frac{1}{\tau(K)}\int_{0}^{\rho_{g_e(K \cap \mathrm{int}\,\mathbb{S}_e^+)}(v)} \frac{r^n}{(1+r^2)^{\frac{n+2}{2}}} \, d{r} = \frac{1}{\tau(L)}\int_{0}^{\rho_{g_e(L \cap \mathrm{int}\,\mathbb{S}_e^+)}(v)} \frac{r^n}{(1+r^2)^{\frac{n+2}{2}}} \, d{r}\]
for all $v\in \mathbb{S}^{n-1}$. Thus, since $t\to \int_{0}^{t} r^n(1+r^2)^{-\frac{n+2}{2}}\, d{r}$ is strictly increasing, it follows that $\rho_{g_e(K \cap \mathrm{int}\,\mathbb{S}_e^+)}(v) \leq \rho_{g_e(L \cap \mathrm{int}\,\mathbb{S}_e^+)}(v)$ for all $v\in \mathbb{S}^{n-1}$ or, equivalently, $K \subseteq L$. Hence, if $K$ and $L$ have equal spherical volume, they must coincide.

\pagebreak

Finally, for the proof of (d) assume that $K \in \mathcal{K}_c(\mathbb{S}^n)$ has center $e \in \mathbb{S}^n$. Since the restriction of $g_e$ to any spherical cap of radius $\alpha < \frac{\pi}{2}$ is a diffeomorphism onto some Euclidean ball in $\mathbb{R}_e^n$, the boundary of $\Gamma_s K$ is a $C^2$ submanifold by Lemma \ref{2:lemMuCentBodC2} (and its extension to unbounded convex sets discussed in the proof of Proposition~\ref{3:scbGnom}). Moreover, it follows from \textbf{\cite[\textnormal{Lemma 4.4}]{Besau:Werner:2016}} that the spherical Gau\ss\ curvature of $\Gamma_s K$ at $u \in \mathbb{S}^n$ vanishes precisely when the one of $\Gamma_{\widehat{\tau}\,} g_e(K \cap \mathrm{int}\, \mathbb{S}_e^+)$ vanishes at $g_e(u) \in \mathbb{R}^n_e$. Hence, Lemma \ref{2:lemMuCentBodC2} and its extension complete the proof. \hfill $\blacksquare$

\vspace{0.3cm}

Before we continue, we remark that, like Blaschke's classical centroid body map (cf.\ \textbf{\cite{Lutwak:1990}}), it is not difficult to see that $\Gamma_s$ is \emph{not} monotone under set inclusion.

\vspace{0.3cm}

In the last part of this section, we establish a couple of auxiliary results concerning the discrete spherical centroid bodies defined in the introduction,
\[\Gamma_{s,e}(u_1,\ldots, u_N) = \mathrm{conv}\left\{c_s\left(u_1^{(e)},\ldots,u_N^{(e)}\right)\right\},\]
which are used in Theorem \ref{1:thmApprox} to approximate spherical centroid bodies of convex bodies. (Recall that here, $\{u_1, \ldots, u_N\} \subseteq \mathrm{int}\,\mathbb{S}_e^+$ for some $e \in \mathbb{S}^n$ and $u^{(e)} = \{u, u^e\}$.)

Note that, by definition and Lemma \ref{2:propSpherCent}, $\Gamma_{s,e}(u_1,\ldots, u_N) \in \mathcal{K}_c(\mathbb{S}^n)$ has center $e$ and is proper. Lemma \ref{2:propSpherCent} also implies that the map $\Gamma_{s,e}$ is continuous and $\mathrm{O}(n+1)$-equivariant. Moreover, as a direct consequence of Lemma \ref{2:propGnom} (d), Lemma \ref{2:csGnom}, and Lemma \ref{2:fCentroidBodies} we obtain for its gnomonic projection the following.

\begin{koro}\label{4:gnomApprox}
For $e \in \mathbb{S}^n$, let $g_e: \mathrm{int}\,\mathbb{S}_e^+ \to \mathbb{R}^n_e$ denote the gnomonic projection and $\{u_1,\ldots, u_N\} \subseteq \mathrm{int}\,\mathbb{S}_e^+$. Then
\[g_e(\Gamma_{s,e}(u_1,\ldots, u_N)) = \Gamma_\phi(g(u_1),\ldots, g(u_N)).\]
\end{koro}

Recall that our definition of $\Gamma_{s,e}(u_1,\ldots, u_N)$ was motivated by relation (\ref{randcent}) for discrete centroid bodies in a linear vector space. However,
in the linear setting, there is an alternative way to express these $\Gamma(x_1,\ldots,x_N)$, for $x_1,\ldots,x_n\in \mathbb{R}^n$, namely,
\[\Gamma(x_1,\ldots, x_N) = \frac{1}{N}\sum_{i=1}^N [-x_i, x_i] = \{c(y_1,\ldots, y_N): y_i\in[-x_i, x_i], 1\leq i\leq N\}.\]
By mimicking this approach on the sphere, we define, $\{u_1,\ldots, u_N\} \subseteq \mathrm{int}\,\mathbb{S}_e^+$,
\[  \widetilde{\Gamma}_{s,e}(u_1,\ldots, u_N) := \{c_s(v_1,\ldots, v_N): v_i\in[u_i^e, u_i], 1\leq i\leq N\},\]
where $[u_i^e, u_i]$ denotes the geodesic segment connecting $u_i^e$ and the geodesic reflection of $u_i$ about $e$. These new sets are, in general, not spherically convex. However, there is the following interesting relation between them and $\Gamma_{s,e}(u_1,\ldots, u_N)$.

\begin{prop} \label{discrete1and2}
For $e \in \mathbb{S}^n$ and any \{$u_1,\ldots,u_N\} \subseteq \mathrm{int}\,\mathbb{S}_e^+$, we have
\[  \Gamma_{s,e}(u_1,\ldots, u_N) = \mathrm{conv}\, \widetilde{\Gamma}_{s,e}(u_1,\ldots, u_N).\]
\end{prop}

\pagebreak

\noindent {\it Proof.} Let $g_e: \mathrm{int}\,\mathbb{S}_e^+ \to \mathbb{R}_e^n$ denote gnomonic projection. Then, by Lemma \ref{2:propGnom}, it suffices to prove that
\[g_e(\Gamma_{s,e}(u_1,\ldots, u_N))= \mathrm{conv}\,g_e(\widetilde{\Gamma}_{s,e}(u_1,\ldots, u_N)). \]
But, by Lemma \ref{4:gnomApprox} and Lemma \ref{2:csGnom}, this is equivalent to
\begin{equation}\label{5:eqApproxEqualGnom}
  \Gamma_\phi(x_1,\ldots,x_N) = \mathrm{conv}\left\{c_\phi(y_1,\ldots, y_N): y_i\in[-x_i, x_i], 1\leq i\leq N\right\},
\end{equation}
where $x_i = g_e(u_i)$, $1 \leq i \leq N$. In order to prove (\ref{5:eqApproxEqualGnom}), note that, by Lemma \ref{2:fCentroidBodies},
\begin{align*}
 \Gamma_\phi(x_1,\ldots, x_N) &= \mathrm{conv}\left\{c_\phi(\pm x_1,\ldots, \pm x_N)\right\} \\
	                             &\subseteq \mathrm{conv} \left\{c_\phi(y_1,\ldots, y_N):  y_i\in[-x_i, x_i], 1\leq i\leq N\right\}.
\end{align*}
Thus, it only remains to prove the reverse inclusion. To this end, recall that for $z_1,\ldots, z_N\in\mathbb{R}^n$ and $v \in \mathbb{S}^n$,
\[  h(\Gamma_\phi(z_1,\ldots, z_N), v) = \frac{1}{\sum_{i=1}^N \phi(z_i)} \sum_{i=1}^N \phi(z_i)|v\cdot z_i|.\]
Using  $\nabla \phi(x) = -\phi(x)^3 x$, a straightforward computation yields
\[ \left . \frac{d}{d\lambda}\right|_{\lambda=1} h(\Gamma_\phi(\lambda z_1,\ldots, z_N), v) =
  \frac{\|z_1\|^2 \phi(z_1)^3}{\left( \sum_{i=1}^N \phi(z_i) \right) ^2} \sum_{i=1}^N \phi(z_i)\left[\frac{|v \cdot z_1|}{\|z_1\|^2}+ |v \cdot z_i| \right] > 0.\]
Repeating this computation for $z_2,\ldots, z_N$ shows that for any $v \in \mathbb{S}^{n-1}$, the function $(z_1,\ldots, z_N)\mapsto h(\Gamma_\phi(z_1,\ldots, z_N), \theta)$ is radially increasing in every coordinate. By applying this fact to each coordinate successively, we obtain for all $y_i\in[0, x_i]$, $1\leq i\leq N$, and every $v \in \mathbb{S}^{n-1}$,
\[h(\Gamma_\phi(y_1,\ldots, y_N), v) \leq h(\Gamma_\phi(x_1,\ldots,x_N), v),\]
that is, $\Gamma_\phi(y_1,\ldots,y_N) \subseteq \Gamma_\phi(x_1,\ldots,x_N)$. But, since both sets are origin-symmetric, this inclusion also holds for all $y_i\in[-x_i, x_i]$, $1\leq i\leq N$. In particular,
\[c_\phi(y_1,\ldots, y_N) \in \Gamma_\phi(y_1,\ldots,y_N) \subseteq \Gamma_\phi(x_1,\ldots,x_N).\]
Since $\Gamma_\phi(x_1,\ldots,x_N)$ is convex, this concludes the proof. \hfill $\blacksquare$

\vspace{0.3cm}

In Section 5, we present the proof of Theorem \ref{1:thmApprox}, showing that the discrete centroid bodies $\Gamma_{s,e}(u_1,\ldots, u_N)$ approximate $\Gamma_s K$, when $u_1,\ldots, u_N$ are chosen randomly from $K$. By Proposition \ref{discrete1and2}, the same holds true for the bodies $\mathrm{conv}\, \widetilde{\Gamma}_{s,e}(u_1,\ldots, u_N)$. Our final result of this section is a critical ingredient in the proof of these facts and based on a variant of the proof of \textbf{\cite[\textnormal{Corollary 5.2}]{Paouris:Pivovarov:2012}}.

\begin{lem}\label{4:approxGammaMu}
Let $\mu, \nu$ be finite, absolutely continuous Borel measures on $\mathbb{R}^n$ and let $f$ denote the density of $\mu$ with respect to $\nu$. Then, for $L \in \mathcal{K}(\mathbb{R}^n)$ and independent random vectors $X_1, \ldots, X_N$ on $\mathbb{R}^n$, identically distributed according to $\frac{\mathbbm{1}_L}{\nu(L)}\,d\nu$, we have
\[\Gamma_f(X_1,\ldots, X_N) \to \Gamma_\mu L\]
almost surely in the Hausdorff metric as $N \rightarrow \infty$.
\end{lem}

\noindent {\it Proof.} By the strong law of large numbers (see e.g., \textbf{\cite[\textnormal{Theorem 8.3.5}]{dudley2002}}), we have
\[\frac{1}{N}\sum_{i=1}^N f(X_i) \to \frac{1}{\nu(L)}\int_L f(x) \, d\nu(x) = \frac{\mu(L)}{\nu(L)}\]
and
\[\frac{1}{N}\sum_{i=1}^N |y\cdot X_i|f(X_i) \to \frac{1}{\nu(L)} \int_L |y\cdot x|f(x) \,d\nu(x) = \frac{1}{\nu(L)} \int_L |y\cdot x| \,d\mu(x)  \]
almost surely for every $y \in \mathbb{R}^n$ as $N \to \infty$. Since the product of almost surely convergent sequences of random variables converges almost surely to the product of their respective limits, we conclude
\begin{align*}
    h(\Gamma_f(X_1,\ldots, X_N), y) & = \frac{1}{\sum_{i=1}^N f(X_i)} \sum_{i=1}^N |y\cdot X_i|f(X_i) \\
                                    &\to \frac{1}{\mu(L)} \int_L |y\cdot x|\, d\mu(x) = h(\Gamma_\mu(L), y),
\end{align*}
almost surely for every $y \in \mathbb{R}^n$. This proves the desired statement, since pointwise convergence of support functions is equivalent to the convergence of the respective bodies in the Hausdorff metric (see e.g., \textbf{\cite[\textnormal{p.\ 54}]{Schneider:2014}}). \hfill $\blacksquare$

\vspace{0.3cm}

Finally, we note that Lemma \ref{4:approxGammaMu} holds true for any closed and unbounded convex set $L$ in $\mathbb{R}^n$ as long as $\mu$ has finite first moments (so that $\Gamma_\mu L$ exists).

\vspace{1cm}

\centerline{\large{\bf{ \setcounter{abschnitt}{4}
\arabic{abschnitt}. Further auxiliary results}}}

\reseteqn \alpheqn \setcounter{theorem}{0}

\vspace{0.6cm}

In this shorter section, we recall additional notions and results required in the proof of Theorem \ref{1:thmIneq}. In particular,  we discuss a crucial result from \textbf{\cite{Cordero:etal:2015}} about the expected behavior of certain random functionals under symmetric decreasing rearrangements.

To this end, first recall that a function $f: \mathbb{R}^n \to \mathbb{R}^+$ is called \emph{radially symmetric}, if $f(x) = \bar{f}(\|x\|)$ for all $x\in \mathbb{R}^n$ for some $\bar{f}: \mathbb{R}^+\to \mathbb{R}^+$. The function $f$ is called \emph{radially decreasing}, if for every $v \in \mathbb{S}^{n-1}$, the functions $f_v: \mathbb{R}^+ \to \mathbb{R}^+$, $f_v(r)=f(rv)$, are decreasing.

\pagebreak

For a bounded Borel set $A \subseteq \mathbb{R}^n$ we write in the following $A^{\mbox{\begin{tiny}$\bigstar$\end{tiny}}}$ to denote the open Euclidean ball centered at the origin such that $\mathrm{vol}(A^{\mbox{\begin{tiny}$\bigstar$\end{tiny}}})=\mathrm{vol}(A)$. Then, for an integrable function $f: \mathbb{R}^n \to \mathbb{R}^+$, its \emph{symmetric decreasing rearrangement} $f^{\mbox{\begin{tiny}$\bigstar$\end{tiny}}} : \mathbb{R}^n \to \mathbb{R}^+$ is defined by
\[f^{\mbox{\begin{tiny}$\bigstar$\end{tiny}}}(x)=\int_0^\infty \mathbbm{1}_{\{f>t\}^{\mbox{\begin{tiny}$\bigstar$\end{tiny}}}}(x)\, dt. \]
Note that the function $f^{\mbox{\begin{tiny}$\bigstar$\end{tiny}}}$ is radially symmetric, radially decreasing, and satisfies
\begin{equation}\label{4:rearrInt}
  \int_{\R^n}f(x)\, dx=\int_{\R^n}f^{\mbox{\begin{tiny}$\bigstar$\end{tiny}}}(x)\, dx.
\end{equation}
Hence, the symmetric decreasing rearrangement of a probability density is again a probability density. For a random vector $X$ in $\mathbb{R}^n$ distributed according to a probability density $f$, we therefore write $X^{\mbox{\begin{tiny}$\bigstar$\end{tiny}}}$ for the random vector in $\mathbb{R}^n$ distributed according to $f^{\mbox{\begin{tiny}$\bigstar$\end{tiny}}}$.

The following theorem is the special case of \textbf{\cite[\textnormal{Theorem 4.2}]{Cordero:etal:2015}}, where we take $r=0$ and $C = \frac{1}{N}B_\infty^N$, the rescaled $\ell_{\infty}^N$ unit ball.

\begin{theorem}[Cordero-Erausquin et al.\ \textbf{\cite[\textnormal{Theorem 4.2}]{Cordero:etal:2015}}]\label{4:symRearr}
Let $\mu, \nu$ be finite, absolutely continuous Borel measures on $\mathbb{R}^n$ such that the density of $\nu$ is radially symmetric and radially decreasing. Then, for $L \in \mathcal{K}(\mathbb{R}^n)$ and independent random vectors $X_1,\ldots, X_N$ in $\mathbb{R}^n$, identically distributed according to $\frac{\mathbbm{1}_L}{\mu(L)}\,d\mu$, we have
\[\mathbb{E}\left[\nu\left(\Gamma^\circ (X_1, \ldots, X_N)\right)\right] \leq \mathbb{E}\left[\nu\left(\Gamma^\circ (X_1^{\mbox{\begin{tiny}$\bigstar$\end{tiny}}}, \ldots, X_N^{\mbox{\begin{tiny}$\bigstar$\end{tiny}}})\right)\right].\]
\end{theorem}

Recall that $\Gamma^\circ (X_1, \ldots, X_N)$ denotes the polar of the origin-symmetric convex body $\Gamma(X_1, \ldots, X_N)$. We also note that the proof of Theorem \ref{4:symRearr} from \textbf{\cite{Cordero:etal:2015}} was given for densities bounded by one. However, it goes through verbatim in the case of arbitrary bounded probability densities.

Next, following ideas from \textbf{\cite{Cordero:etal:2015}}, we apply a bathtub-type argument to transform the rearrangement inequality from Theorem \ref{4:symRearr} into a randomized isoperimetric inequality. To this end, recall that for $L \in \mathcal{K}(\mathbb{R}^n)$ and a finite Borel measure $\mu$ on $\mathbb{R}^n$, we write $B_L^{\mu}$ for the Euclidean ball centered at the origin such that $\mu(L)=\mu(B_L^\mu)$. In addition, we use in the following $Z_1, \ldots, Z_N$ to denote independent random vectors in $\mathbb{R}^n$ which are identically distributed according to \[\frac{\mathbbm{1}_{B_L^\mu}}{\mu(B_L^\mu)}\, d\mu.\]

\begin{theorem}\label{4:massTrans}
Let $\mu, \nu$ be finite, absolutely continuous Borel measures on $\mathbb{R}^n$ such that the density of $\mu$ is radially symmetric and radially decreasing. Then, for \linebreak $L \in \mathcal{K}(\mathbb{R}^n)$ and independent random vectors $X_1,\ldots, X_N$ on $\mathbb{R}^n$, identically distributed according to $\frac{\mathbbm{1}_L}{\mu(L)}\,d\mu$, we have
\[\mathbb{E}\left[\nu\left(\Gamma^\circ (X_1^{\mbox{\begin{tiny}$\bigstar$\end{tiny}}}, \ldots, X_N^{\mbox{\begin{tiny}$\bigstar$\end{tiny}}})\right)\right] \leq	\mathbb{E}\left[\nu\left(\Gamma^\circ (Z_1, \ldots, Z_N)\right)\right].\]
\end{theorem}

\noindent {\it Proof.} If $f_\mu$ denotes the density of $\mu$, then, by (\ref{4:rearrInt}), we have
\[ \int_{\R^n}(\mathbbm{1}_{B_L^\mu} f_\mu)(x)\, dx = \mu(B_L^\mu) = \mu(L) = \int_{\R^n}(\mathbbm{1}_L f_\mu)(x)\, dx = \int_{\R^n}(\mathbbm{1}_L f_\mu)^{\mbox{\begin{tiny}$\bigstar$\end{tiny}}}(x)\, dx.\]
Since $f_\mu$ is radially symmetric, there exists $\overline{f_{\mu}}: \mathbb{R}^+ \to \mathbb{R}^+$ such that $\overline{f_\mu}(\|x\|) = f_\mu(x)$ for all $x \in \mathbb{R}^n$.
Hence, denoting by $r_L^\mu$ the radius of the ball $B_L^\mu$ and using polar coordinates, we see that
\begin{equation} \label{masstrans17}
\int_0^{r_L^\mu} \overline{f_\mu}(r) r^{n-1}\, dr = \int_0^\infty \overline{(\mathbbm{1}_L f_\mu)^{\mbox{\begin{tiny}$\bigstar$\end{tiny}}}}(r) r^{n-1}\, dr.
\end{equation}
Now, define a function $\alpha:\mathbb{R}^+\to\mathbb{R}$ by
\[ \alpha(r):=\left (\mathbbm{1}_{[0,r_L^\mu]}(r)\overline{f_\mu}(r)-\overline{(\mathbbm{1}_Lf_\mu)^{\mbox{\begin{tiny}$\bigstar$\end{tiny}}}}(r)\right)r^{n-1}.\]
Then, by (\ref{masstrans17}) and the fact that $(\mathbbm{1}_L f_\mu)^{\mbox{\begin{tiny}$\bigstar$\end{tiny}}}\leq f_\mu$ (because $f_\mu$ is radially symmetric and radially decreasing), $\alpha$ has the following two properties
\[\mbox{(i) } \int_0^\infty \alpha(r)\, dr=0 \quad \qquad \quad \mbox{(ii)  } \alpha(r) \left \{ \begin{array}{l} \leq 0 \mbox{ for } r > r_L^\mu, \\ \geq 0 \mbox{ for } r \leq r_L^\mu.  \end{array}  \right .  \]
Combining (i) and (ii), it follows that for any radially decreasing $F:\mathbb{R}^n\to\mathbb{R}^+$, we have
\[\int_0^\infty \overline{F}(r)\alpha(r)\, dr = \int_0^\infty (\overline{F}(r)-\overline{F}(r_L^\mu))\alpha(r)\, dr \geq 0\]
or, equivalently by the definition of $\alpha$, that
\[\int_0^\infty \overline{F}(r)\overline{(\mathbbm{1}_L f_\mu)^{\mbox{\begin{tiny}$\bigstar$\end{tiny}}}}(r)r^{n-1}\, dr \leq \int_0^{r_L^\mu} \overline{F}(r)\overline{f_\mu}(r)r^{n-1}\, dr.\]
Transferring back to cartesian coordinates, this inequality becomes
\[\int_{\R^n} F(x)(\mathbbm{1}_L(x)f_\mu(x))^{\mbox{\begin{tiny}$\bigstar$\end{tiny}}}\, dx\leq \int_{\R^n} F(x)\mathbbm{1}_{B_L^\mu}(x)f_\mu(x)\, dx.\]
Now, given $F: (\mathbb{R}^n)^N\to\mathbb{R}^+$ that is radially decreasing in each coordinate, we can apply the above inequality coordinatewise and use Fubini's theorem to obtain
\begin{align*}
  \int_{(\R^n)^N} &F(x_1,\ldots, x_N)
    \prod_{i=1}^N \left(\frac{\mathbbm{1}_L(x_i)}{\mu(L)}f_\mu(x_i)\right)^{\mbox{\begin{tiny}$\bigstar$\end{tiny}}}\, dx_1\cdots dx_N \\
	                &\leq \int_{(\R^n)^N} F(x_1,\ldots, x_N) \prod_{i=1}^N \frac{\mathbbm{1}_{B_L^\mu}(x_i)}{\mu(L)}f_\mu(x_i)\, dx_1\cdots dx_N.
\end{align*}
Finally, putting $F(x_1,\ldots, x_N) = \nu(\Gamma^\circ(x_1,\ldots, x_N))$ yields the desired inequality. \hfill $\blacksquare$

\pagebreak

\centerline{\large{\bf{ \setcounter{abschnitt}{5}
\arabic{abschnitt}. Proofs of the main results}}}

\reseteqn \alpheqn \setcounter{theorem}{0}

\vspace{0.6cm}

The careful preparations of the last two sections allow us now to complete the proofs of Theorems \ref{1:thmApprox} and \ref{1:thmIneq} almost effortlessly. The end of this
final section contains a couple of concluding remarks and natural open problems concerning spherical centroid bodies.

We begin with the proof of Theorem \ref{1:thmApprox} and also recall the statement for the reader's convenience.

\begin{theorem} \label{thm1ext}
Let $K \in \mathcal{K}_c(\mathbb{S}^n)$ have center $e \in \mathbb{S}^n$. If $U_1, \ldots, U_N$ are independent random unit vectors distributed uniformly in $K$, then
\[\Gamma_{s,e}(U_1,\ldots,U_N) \to \Gamma_s K\]
almost surely in the spherical Hausdorff metric as $N \to \infty$.
\end{theorem}

\noindent {\it Proof.} Let us first assume that $K$ is proper, that is, $K \subseteq \mathrm{int}\,\mathbb{S}_e^+$, and let again $g_e : \mathrm{int}\,\mathbb{S}_e^+ \to \mathbb{R}_e^n$ denote the gnomonic projection. Putting $X_i := g_e(U_i)$, $1\leq i\leq N$, the independence and uniform distribution of the $U_i$ together with Lemma \ref{2:pushMeasGnom} (a), implies that the $X_i$ are independent random vectors in $\mathbb{R}^n_e$, identically distributed according to
\[\frac{\mathbbm{1}_{g_e(K)}}{\widehat{\sigma}(g_e(K))}\, d\phantom{.}\!\widehat{\sigma}.\]
Moreover, by Lemma \ref{2:pushMeasGnom}, we have
\[\phi(x) = (1 + \|x\|^2)^{-\frac{1}{2}} = \frac{(1 + \|x\|^2)^{-\frac{n+2}{2}}}{(1 + \|x\|^2)^{-\frac{n+1}{2}}}= \frac{d\phantom{.}\!\widehat{\tau}}{d\phantom{.}\!\widehat{\sigma}}.\]
Hence, by Lemma \ref{4:approxGammaMu}, we obtain
\[\Gamma_\phi(X_1,\ldots, X_N) \to \Gamma_{\widehat{\tau}}\left(g_e(K)\right)\]
almost surely in the Hausdorff metric as $N \to \infty$. Applying now $g_e^{-1}$ and using Lemma \ref{2:propGnom} and Corollaries \ref{3:scbGnom2} and \ref{4:gnomApprox}, we arrive at the desired statement
\[\Gamma_{s,e}(U_1,\ldots, U_N) \to \Gamma_s K.\]

If $K$ is not proper, then we still have $K \subseteq \mathbb{S}_e^+$, and, since $K \setminus\mathrm{int}\,\mathbb{S}_e^+$ is a set of measure zero, we may assume that $U_1,\ldots, U_N$ lie in
$\mathrm{int}\,\mathbb{S}_e^+$. Thus, as in the first part of the proof, it follows from Lemma \ref{4:approxGammaMu} and the remark directly following it, that
\[\Gamma_\phi(X_1,\ldots, X_N) \to \Gamma_{\widehat{\tau}\!}\left(g_e(K \cap \mathrm{int}\,\mathbb{S}_e^+)\right),\]
where, as before, $X_i := g_e(U_i)$, $1\leq i\leq N$. Applying $g_e^{-1}$ and using Corollaries \ref{3:scbGnom2} and \ref{4:gnomApprox} yields again the desired result.
\hfill $\blacksquare$

\vspace{0.3cm}

As an immediate consequence of Theorem \ref{thm1ext} and Proposition \ref{discrete1and2}, we note the following:

\begin{koro}
Let $K \in \mathcal{K}_c(\mathbb{S}^n)$ have center $e \in \mathbb{S}^n$. If $U_1, \ldots, U_N$ are independent random unit vectors distributed uniformly in $K$, then
\[\mathrm{conv}\,\widetilde{\Gamma}_{s,e}(U_1,\ldots, U_N) \to \Gamma_s K\]
almost surely in the spherical Hausdorff metric as $N \to \infty$.
\end{koro}

We turn to the proof of Theorem \ref{1:thmIneq} which is based on the following proposition of independent interest.

\begin{prop}\label{5:lemIneqEucl}
Let $\mu, \nu$ be finite, absolutely continuous Borel measures on $\mathbb{R}^n$ such that their densities are radially symmetric and radially decreasing. Then, for an origin-symmetric convex body $L \in \mathcal{K}(\mathbb{R}^n)$, we have
\begin{equation*}
\nu(\Gamma_\mu^\circ L) \leq \nu(\Gamma_\mu^\circ B_L^\mu).
\end{equation*}
\end{prop}

\noindent {\it Proof.} Let $X_1,\ldots, X_N$ and $Z_1,\ldots, Z_N$ be two families of independent random vectors in $\mathbb{R}^n$ such that each family is identically distributed according to
\begin{align*}
 \frac{\mathbbm{1}_L(x)}{\mu(L)} d\mu(x) \qquad \mbox{ and } \qquad  \frac{\mathbbm{1}_{B_L^\mu}(x)}{\mu(B_L^\mu)} d\mu(x),
\end{align*}
respectively. Then, combining Theorems \ref{4:symRearr} and \ref{4:massTrans}, we obtain
\begin{equation}\label{5:eqIneqN}
  \E\left[\nu\left(\Gamma^\circ(X_1, \ldots, X_N)\right)\right] \leq
	\E\left[\nu\left(\Gamma^\circ(Z_1\ldots Z_N)\right)\right].
\end{equation}
Now, by Proposition \ref{4:approxGammaMu}, we know that $\Gamma(X_1, \ldots, X_N) \to \Gamma_\mu(L)$ almost surely in the Hausdorff metric as $N\to\infty$. Moreover, since taking the polar body and $\nu$ are continuous on origin-symmetric convex bodies in $\mathbb{R}^n$ (see \textbf{\cite[\textnormal{Lemma 5.2}]{Cordero:etal:2015}}), we also have that
\[\nu\left(\Gamma^\circ(X_1, \ldots, X_N)\right) \to \nu\left(\Gamma_\mu^\circ L\right) \]
almost surely as $N\to\infty$. Since $\Gamma_\mu L$ has nonempty interior, there exists $r>0$ such that, for $N$ large enough, we have $rB_2^n \subseteq \Gamma(X_1, \ldots, X_N)$ almost surely and, hence, $\nu\left(\Gamma^\circ(X_1, \ldots, X_N)\right)\leq \nu\left(\frac{1}{r}B_2^n\right)$ almost surely. Therefore, by the theorem of dominated convergence, we conclude that
\[ \E\left[\nu\left(\Gamma^\circ(X_1, \ldots, X_N)\right)\right] \to \E\left[\nu\left(\Gamma_\mu^\circ L \right)\right] = \nu\left(\Gamma_\mu^\circ L \right)\]
and, by the same arguments,
\[\E\left[\nu\left(\Gamma^\circ(Z_1, \ldots, Z_N)\right)\right] \to \E\left[\nu\left(\Gamma_\mu^\circ B_L^\mu\right)\right] = \nu\left(\Gamma_\mu^\circ B_L^\mu\right).\]
Thus, by letting $N \to \infty$ in (\ref{5:eqIneqN}), we obtain the desired inequality. \hfill $\blacksquare$

\vspace{0.3cm}

We are now in a position to prove Theorem \ref{1:thmIneq}:

\begin{theorem}
If $K \in \mathcal{K}_c(\mathbb{S}^n)$ has center $e \in \mathbb{S}^n$, then
\[\sigma(\Gamma_s^* K) \leq \sigma(\Gamma_s^* C_K^\tau).\]
\end{theorem}

\pagebreak

\noindent {\it Proof.}
Let us first assume that $K$ is proper, that is, $K \subseteq \mathrm{int}\,\mathbb{S}_e^+$, and let again $g_e: \mathrm{int}\,\mathbb{S}_e^+ \to \mathbb{R}_e^n$ denote gnomonic projection. Since, by Lemma \ref{2:pushMeasGnom}, the push-forwards $g_e\#\sigma=:\widehat{\sigma}$ and $g_e\#\tau=:\widehat{\tau}$ have radially symmetric and radially decreasing densities,
an application of Proposition \ref{5:lemIneqEucl} to the origin-symmetric convex body $g_e(K)$ yields
\begin{equation}\label{5:eqIneqTHM1EuclApplied}
\widehat{\sigma}\left( \Gamma_{\widehat{\tau}\,}^\circ g_e(K) \right) \leq \widehat{\sigma}\left( \Gamma_{\widehat{\tau}\,}^\circ B_{g_e(K)}^{\widehat{\tau}} \right)=
\widehat{\sigma}\left( \Gamma_{\widehat{\tau}\,}^\circ g_e(C_K^\tau) \right).
\end{equation}
Thus, using Corollary \ref{3:scbGnom2}, the desired inequality
\[\sigma(\Gamma_s^* K) \leq \sigma(\Gamma_s^* C_K^\tau)\]
follows by applying $g_e^{-1}$ to (\ref{5:eqIneqTHM1EuclApplied}).

If $K \in \mathcal{K}_c(\mathbb{S}^n)$ is nonproper, then we still have $K \subseteq \mathbb{S}_e^+$, and we can choose a sequence $K_m \in \mathcal{K}_c(\mathbb{S}^n)$ of proper convex bodies with center $e$ such that $K_m \to K$ in the spherical Hausdorff metric. Moreover, by the first part of the proof, we know that $\sigma(\Gamma_s^* K_m) \leq \sigma(\Gamma_s^* C_{K_m}^\tau)$ for all $m\in\N$. But, since $\sigma$, $\tau$, and $\Gamma_s$ are continuous on $\mathcal{K}_c(\mathbb{S}^n)$ (the latter by Proposition \ref{3:scbProp} (b)), and taking the polar is continuous on proper bodies in $\mathcal{K}_c(\mathbb{S}^n)$ (recall that $\Gamma_sL$ is proper for all $L \in \mathcal{K}_c(\mathbb{S}^n)$ by Corollary~\ref{3:scbGnom2}), we obtain the desired inequality by letting $m \to \infty$.
\hfill $\blacksquare$

\vspace{0.3cm}

We conclude the article with three remarks concerning possible extensions and improvements of Theorem \ref{1:thmIneq}. We begin by discussing a version for not necessarily centrally-symmetric bodies. To this end let $K \in \mathcal{K}(\mathbb{S}^n)$ be proper and assume that $K \subseteq \mathrm{int}\, \mathbb{S}_w^+$ for some $w \in \mathbb{S}^n$ or, equivalently, that $w \in - \mathrm{int}\,K^*$. Moreover, let $\tau_w$ denote the absolutely continuous measure on $\mathbb{S}^n$ with density $u \mapsto |u\cdot w|$ and let $\widehat{\tau}_w := g_w \# \tau_w$ denote its pushforward under gnomonic projection $g_w: \mathrm{int}\,S_w^+ \rightarrow \mathbb{R}^n_w$. If we define the spherical centroid body of $K$ by
\begin{equation} \label{asymmdef}
g_w(\Gamma_s K) = \Gamma_{\widehat{\tau}_w} g_w(K),
\end{equation}
then the arguments leading up to Theorem \ref{1:thmIneq} yield the inequality
\begin{align*}
\sigma(\Gamma_s^* K) \leq \sigma(\Gamma_s^* C_K^{\tau_w}).
\end{align*}
However, we are reluctant to use (\ref{asymmdef}) as definition for $\Gamma_s K$, since, on the one hand, it is not intrinsic and, on the other hand, there is the question what would be a natural choice for $w\in -\mathrm{int}\, K^*$? Of course, this choice should coincide with the center for centrally-symmetric bodies, like, for example, the centroid $c_s(K)$. But we do not know whether $c_s(K) \in -\mathrm{int}\, K^*$ for every proper $K \in \mathcal{K}(\mathbb{S}^n)$, when $n \geq 3$.

\vspace{0.2cm}

Our second remark concerns a possible version of Theorem \ref{1:thmIneq}, where $C_K^\tau$ is replaced by $C_K^\sigma$, that is, the inequality
\begin{equation} \label{cksigma}
\sigma(\Gamma_s^* K) \leq \sigma(\Gamma_s^* C_K^\sigma).
\end{equation}
This would be a stronger isoperimetric inequality than Theorem \ref{1:thmIneq}, since $C_K^{\tau} \subseteq C_K^\sigma$ for every $K \in \mathcal{K}_c(\mathbb{S}^n)$ with center $e \in \mathbb{S}^n$ and equality holds if and only if $K$ is already a cap centered at $e$.

\pagebreak

A possible approach to establishing (\ref{cksigma}) is via a spherical analogue of inequality (\ref{5:eqIneqN}). More precisely, if $U_1, \ldots, U_N$ and $V_1, \ldots, V_N$ are independent random unit vectors uniformly distributed in $K$ and $C_K^{\sigma}$, respectively, is it true that
\begin{equation} \label{randsphpolBPCinequ}
\E[\sigma\left(\Gamma_{s,e}^*(U_1,\ldots, U_N)\right)] \leq \E[\sigma\left(\Gamma_{s,e}^*(V_1,\ldots, V_N)]\right)?
\end{equation}
A combination of inequality (\ref{randsphpolBPCinequ}) with Theorem \ref{1:thmApprox} would then yield (\ref{cksigma}).

\vspace{0.2cm}

Finally, let us state the most interesting and probably hardest open problem concerning spherical centroid bodies -- a spherical analogue of the Busemann--Petty centroid inequality:

\vspace{0.25cm}

\noindent {\bf Open Problem.} If $K \in \mathcal{K}_c(\mathbb{S}^n)$, then
\begin{equation} \label{sphBPcentineq}
\sigma(\Gamma_s K) \geq \sigma(\Gamma_s C_K^\sigma).
\end{equation}

\vspace{0.15cm}

Let us emphasize that inequality (\ref{sphBPcentineq}) would not only imply Theorem \ref{1:thmIneq}, by combining (\ref{sphBPcentineq}) with the spherical Blaschke--Santal\'o inequality from \textbf{\cite{Gao:Hug:Schneider:2003}}, but the stronger inequality discussed in the above remark. Moreover, (\ref{sphBPcentineq}) would also imply the classical Busemann--Petty centroid inequality by considering spheres with radii going to infinity and rescaling.

\vspace{0.6cm}

\noindent {{\bf Acknowledgments} T.\ Hack and F.E.\ Schuster were
supported by the European Research Council (ERC), Project number: 306445, and the Austrian Science Fund (FWF), Project numbers:
Y603-N26 and P31448-N35. F.\ Besau was partially supported by the Deutsche Forschungsgemeinschaft (DFG), Project number: BE2484/5-2.
P.\ Pivovarov was supported by the NSF grant DMS-1612936.

\vspace{-0.3cm}

\begin{small}

\[ \begin{array}{ll}
\mbox{Florian Besau}                                     & \mbox{Thomas Hack} \\
\mbox{Vienna University of Technology \phantom{wwwwWWW}} & \mbox{Vienna University of Technology} \\
\mbox{florian.besau@tuwien.ac.at}                        & \mbox{thomas.hack@tuwien.ac.at} \\
& \\
\mbox{Peter Pivovarov}                                   & \mbox{Franz Schuster} \\
\mbox{University of Missouri \phantom{wwwwWW}}           & \mbox{Vienna University of Technology} \\
\mbox{pivovarovp@missouri.edu}                           & \mbox{franz.schuster@tuwien.ac.at}
\end{array}\]

\end{small}


\begin{thebibliography}{99}
\footnotesize{
\parskip-0.1cm{

\bibitem{Bernig:2014}
A. Bernig, \emph{Centroid bodies and the convexity of area functionals}, J.\ Differential Geom.\ {\bf 98} (2014), 357--373.
	
\bibitem{Besau:Schuster:2016}
F. Besau and F.E. Schuster, \emph{Binary operations in spherical convex geometry}, Indiana Univ.\ Math.\ J.\ {\bf 65} (2016), 1263--1288.

\bibitem{Besau:Werner:2016}
F. Besau and E.M. Werner, \emph{The spherical convex floating body}, Adv.\ Math.\ {\bf 301} (2016), 867--901.

\bibitem{Besau:Werner:2018}
F. Besau and E.M. Werner, \emph{The floating body in real space forms}, J.\ Differential Geom.\ {\bf 110} (2018), 187--220.

\bibitem{Besau:Ludwig:Werner:2018}
F. Besau, M. Ludwig, and E.M. Werner, \emph{Weighted floating bodies and polytopal approximation}, Trans. Amer. Math. Soc. {\bf 370} (2018), 7129--7148.

\bibitem{Blaschke:1917}
W. Blaschke, \emph{Affine Geometrie IX: Verschiedene Bemerkungen und Aufgaben}, Ber. Verh. S\"achs. Akad. Leipzig. Math.-Phys. Kl. {\bf 69} (1917), 412--420.

\bibitem{Brazitikos:etal:2014}
S. Brazitikos, A. Giannopoulos, P. Valettas, and B.H. Vritsiou, \emph{Geometry of isotropic convex bodies}, Mathematical Surveys and Monographs 196, Amer. Math. Society (2014).

\bibitem{Busemann:1953}
H. Busemann, \emph{Volume in terms of concurrent cross-sections}, Pacific J. Math. {\bf 3} (1953), 1--12.
	
\bibitem{Busemann:1958}
H. Busemann, \emph{Convex surfaces}, Interscience Tracts in Pure and Applied Mathematics, no. 6, Interscience Publishers, Inc., New York; Interscience Publishers Ltd., London, 1958, ix+196.

\bibitem{Campi:Gronchi02a}
S. Campi and P. Gronchi, \emph{The $L_p$-Busemann-Petty centroid inequality}, Adv.\ Math.\ {\bf 167} (2002), 128--141.
	
\bibitem{Cordero:etal:2015}
D. Cordero-Erausquin, M. Fradelizi, G. Paouris, and P. Pivovarov, \emph{Volume of the polar of random sets and shadow systems}, Math.\ Ann.\ {\bf 362} (2015), 1305--1325.

\bibitem{Dann:etal:2018}
S. Dann, J. Kim, and V. Yaskin, \emph{Busemann's intersection inequality in hyperbolic and spherical spaces}, Adv. Math. {\bf 326} (2018), 521--560.

\bibitem{Dann:etal:2016}
S. Dann, G. Paouris, and P. Pivovarov, \emph{Bounding marginal densities via affine isoperimetry}, Proc. Lond. Math. Soc. {\bf 113} (2016), 140--162.

\bibitem{DeNapoli:etal:2018}
P.L. De N\'apoli, J. Haddad, C.H. Jim\'enez, and M. Montenegro, \emph{The sharp affine $L^2$ Sobolev trace inequality and variants}, Math.\ Ann.\ {\bf 370} (2018), 287--308.

\bibitem{dudley2002}
R.M. Dudley, \emph{Real Analysis and Probability}, Cambridge Studies in Advanced Mathematics {\bf 74}, Cambridge University Press, Cambridge, 2002.
	
\bibitem{Galperin:1993}
G.A. Galperin, \emph{A concept of the mass center of a system of material points in the constant curvature spaces}, Comm.\ Math.\ Phys.\ {\bf 154} (1993), 63--84.
	
\bibitem{Gao:Hug:Schneider:2003}
F. Gao, D. Hug, and R. Schneider, \emph{Intrinsic volumes and polar sets in spherical space}, Math.\ Notae\ {\bf 41} (2003), 159--176.
	
\bibitem{Gardner:2006}
R.J. Gardner, \emph{Geometric tomography}, Second edition. Encyclopedia of Mathematics and its Applications 58, Cambridge University Press, New York, 2006.
	
\bibitem{Glasauer:1996}
S. Glasauer, \emph{Integralgeometrie konvexer K\"orper im sph\"arischen Raum}, PhD Thesis, Univ.\ Freiburg, 1995.

\bibitem{Haberl:2012}
C. Haberl, \emph{Minkowski valuations intertwining the special linear group}, J. Eur. Math. Soc. (JEMS) {\bf 14} (2012), 1565--1597.

\bibitem{habschu09}
C. Haberl and F.E. Schuster, \emph{General $L_p$ affine isoperimetric inequalities}, J.\ Differential Geom.\ {\bf 83} (2009), 1--26.

\bibitem{habschu19}
C. Haberl and F.E. Schuster, \emph{Affine vs. Euclidean isoperimetric inequalities}, arXiv:1804.11165.

\bibitem{HJM:2016}
J. Haddad, C.H. Jim\'enez, and M. Montenegro, \emph{Sharp affine Sobolev type inequalities via the $L_p$ Busemann-Petty centroid inequality}, J.\ Funct.\ Anal.\ {\bf 271} (2016), 454--473.

\bibitem{HJM:2018}
J. Haddad, C.H. Jim\'enez, and M. Montenegro, \emph{Sharp affine weighted $L^p$ Sobolev type inequalities}, Trans. Amer. Math. Soc., in press.

\bibitem{HJM:2019}
J. Haddad, C.H. Jim\'enez, and M. Montenegro, \emph{Asymmetric Blaschke-Santal\'o functional inequalities}, arXiv:1810.02288.

\bibitem{Huang:Slomka:Werner:2018}
H. Huang, B.A. Slomka, and E.M. Werner, \emph{Ulam floating bodies}, arXiv:1803.08224.

\bibitem{Ivaki2016}
M.N. Ivaki, \emph{The planar Busemann-Petty centroid inequality and its stability}, Trans. Amer. Math. Soc. {\bf 368} (2016), 3539--3563.

\bibitem{Ivaki2017}
M.N. Ivaki, \emph{The second mixed projection problem and the projection centroid conjectures}, J. Funct. Anal. {\bf 272} (2017), 5144--5161.

\bibitem{Klartag:Milman:2012}
B. Klartag and E. Milman, \emph{Centroid bodies and the logarithmic Laplace transform -- a unified approach}, J.\ Funct.\ Anal.\ {\bf 262} (2012), 10--34.

\bibitem{Ludwig:2007}
M. Ludwig, \emph{Minkowski valuations}, Trans. Amer. Math. Soc. {\bf 357} (2005), 4191--4213.

\bibitem{Lutwak:1986}
E. Lutwak, \emph{Volume of mixed bodies}, Trans. Amer. Math. Soc. {\bf 294} (1986), 487--500.

\bibitem{Lutwak:1990}
E. Lutwak, \emph{Centroid bodies and dual mixed volumes}, Proc. London Math. Soc. {\bf 60} (1990), 365--391.

\bibitem{LZ1997}
E. Lutwak and G. Zhang, \emph{Blaschke-Santal\'o inequalities}, J.\ Differential Geom.\ {\bf 45} (1997), 1--16.

\bibitem{LYZ2000a}
E. Lutwak, D. Yang, and G. Zhang, \emph{$L_p$ affine isoperimetric inequalities}, J.\ Differential Geom.\ {\bf 56} (2000), 111--132.

\bibitem{LYZ2002}
E. Lutwak, D. Yang, and G. Zhang, \emph{Sharp affine $L_p$ Sobolev inequalities}, J.\ Differential Geom.\ {\bf 62} (2002), 17--38.
	
\bibitem{LYZ:2010}
E. Lutwak, D. Yang, and G. Zhang, \emph{Orlicz centroid bodies}, J.\ Differential Geom.\ {\bf 84} (2010), 365--387.
	
\bibitem{Maggi:2012}
F. Maggi, \emph{Sets of Finite Perimeter and Geometric Variational Problems}, An Introduction to Geometric Measure Theory, Cambridge Studies in Advanced Mathematics 135.\ Cambridge University Press, Cambridge, 2012.

\bibitem{Milman:Pajor:1989}
V.D. Milman and A. Pajor, \emph{Isotropic position and inertia ellipsoids and zonoids of the unit ball of a normed $n$-dimensional space}, Geometric Aspects of Functional Analysis (J. Lindenstrauss
and V.D. Milman, eds.) Springer Lecture Notes Math. 1376 (1989), 64--104.
	
\bibitem{Nguyen:2016}
V.H. Nguyen, \emph{New approach to the affine P\'olya-Szeg\"o principle and the stability version of the affine Sobolev inequality}, Adv.\ Math.\ {\bf 302} (2016), 1080--1110.

\bibitem{Nguyen:2018}
V. H. Nguyen, \emph{Orlicz-Lorentz centroid bodies}, Adv.\ in Appl.\ Math.\ {\bf 92} (2018), 99--121.

\bibitem{Oliker:2007}
V. Oliker, \emph{Embedding $\mathbf{S}^n$ into $\mathbf{R}^{n+1}$ with given integral Gauss curvature and optimal mass transport on $\mathbf{S}^n$}, Adv. Math. {\bf 213} (2007), 600--620.

\bibitem{Paouris:2006}
G.$\!$ Paouris, \emph{Concentration of mass on convex bodies}, Geom.$\!$ Funct.$\!$ Anal.$\!$ {\bf 16} (2006), 1021-1049.
	
\bibitem{Paouris:2012}
G. Paouris, \emph{On the existence of supergaussian directions on convex bodies}, Mathematika {\bf 58} (2012), 389--408.

\bibitem{Paouris:2012a}
G. Paouris, \emph{Small ball probability estimates for log-concave measures}, Trans. Amer. Math. Soc. {\bf 364} (2012), 287--308.

\bibitem{Paouris:Pivovarov:2012}
G. Paouris and P. Pivovarov, \emph{A probabilistic take on isoperimetric-type inequalities}, Adv.\ Math.\ {\bf 230} (2012), 1402--1422.

\bibitem{Paouris:Pivovarov:2018}
G. Paouris and P. Pivovarov, \emph{Randomized isoperimetric inequalities}, Convexity and concentration, 391--425, IMA Vol. Math. Appl. 161, Springer, New York, 2017.

\bibitem{Paouris:Valettas:2014}
G. Paouris and P. Valettas, \emph{Neighborhoods on the Grassmannian of marginals with bounded isotropic constant}, J. Funct. Anal. {\bf 267} (2014), 3427--3443.
	
\bibitem{Paouris:Werner:2012}
G. Paouris and E.M. Werner, \emph{Relative entropy of cone measures and $L_p$ centroid bodies}, Proc.\ Lond.\ Math.\ Soc.\ {\bf 104} (2012), 253--286.
	
\bibitem{Petty:1961}
C.M. Petty, \emph{Centroid surfaces}, Pacific J.\ Math.\ {\bf 11} (1961), 1535--1547.
 	
\bibitem{Schneider:2014}
R. Schneider, \emph{Convex bodies: the Brunn-Minkowski theory}, Second expanded edition. Encyclopedia of Mathematics and its Applications 151.\ Cambridge University Press, Cambridge, 2014.

\bibitem{Yaskin:2006}
V. Yaskin, \emph{The Busemann-Petty problem in hyperbolic and spherical spaces}, Adv. Math. {\bf 203} (2006), 537--553.
	
\bibitem{Yaskin:Yaskina:2006}
V. Yaskin and M. Yaskina, \emph{Centroid bodies and comparison of volumes}, Indiana Univ.\ Math.\ J.\ {\bf 55} (2006), 1175--1194.

\bibitem{zhang99}
G. Zhang, \emph{The affine Sobolev inequality}, J.\ Differential Geom.\ {\bf 53} (1999), 183--202.
	
\bibitem{zhu2012}
G. Zhu, \emph{The Orlicz centroid inequality for star bodies}, Adv. Appl. Math. {\bf 48} (2012), 432--445.


}}
\end{thebibliography}
\end{document}